Proofs of the Twin Primes and Goldbach Conjectures
T.J. Hoskins
July 27, 2019

**Abstract**

Associate a unique numerical sequence called the modular signature with each positive integer, using modular residues of each integer under the prime numbers, and distinguishing between the "core seed primes" and "non-core seed primes" used to create the modular signatures. Group the modular signatures within primorials. Use elementary sieve properties and combinatorial principles to prove the twin primes conjecture and the Goldbach conjecture.

## **Introduction**

The twin primes conjecture arose from an open question about the distribution of prime numbers. The conjecture states the following: There are infinitely many integers n, such that n-1 and n+1 are both prime [1]. Stated another way, there are infinitely many pairs of primes whose gap is 2. For example, the twin prime pair 17 and 19 have a gap of 19-17 = 2.

While mathematicians have found many patterns in the prime numbers, they have continued to pursue proofs of the twin primes conjecture. Neale has described the efforts in the last decade to reduce the size of the twin primes gap to 6, conditioned on the truth of the generalized Elliott-Halberstam Conjecture [2].

The Goldbach conjecture is older and arose out of correspondence between Christian Goldbach (1690-1764) and Leonhard Euler (1707-1783) in 1742. In a letter to Euler, Goldbach proposed that every integer greater than 2 could be written as the sum of three primes. Since many mathematicians considered 1 to be a prime number at that time, the modern version of Goldbach's conjecture is that every even integer greater than 4 can be written as the sum of two odd primes [3].

The problem is that if E = 2N is an even number, and a prime number P1 < N, then it is difficult to prove that a prime number P2 exists such that P1 + P2 = 2N for all even numbers. Mathematicians have checked the validity of Goldbach's conjecture using sieve methods. As of 2012, Oliveira e Silva verified that the Goldbach hypothesis is true for N ≤ 4 x $10^{18}$. [4] However, that does not prove the conjecture, because there could be an even number greater than 4 quintillion for which the conjecture is false, given the increasing dispersion of the primes as N increases.

Our approach is narrow. We seek only to prove the existence of an infinite number of twin primes, not to provide a formula to count them. Similarly, we seek only to prove that for any even number, there exists at least one prime pair, P1 and P2, such that P1 + P2 = E, not to provide a formula to count them.
- In Part 1, we will lay the foundation for the proofs by defining the seed primes used in congruence classes to create modular signatures for every positive integer, examining the repetition of modular residues in the modular signatures within primorials, and proving several theorems regarding the dispersion of potential primes within primorials.
- In Part 2, we identify how to predict potential twin primes from modular signatures, examine the dispersion of potential twin primes within the set of potential primes in one or more primorials, and prove theorems on the number of potential twin primes in stacked primorials.
- In Part 3 we develop a rule for identifying Goldbach solutions and describe how to generate potential Goldbach solutions from modular residues, using an equation from Part 2. Additionally, we prove a theorem that the potential primes in a primorial will remain primes if they are less than the square of the smallest non-core seed prime of the primorial.

- In Part 4, we present several theorems that relate to the distribution of potential solutions (twin prime or Goldbach) within a scaffold of primorials and prove two theorems on the existence of Goldbach solutions.
- In Part 5, we present the proof of the twin primes conjecture in which (1) we assume there is a largest prime number, N, that is the last twin prime and (2) create a scaffold of three primorials greater than N, such that the middle primorial is large enough to contain potential twin primes greater than N that are less than the square of the smallest non-core seed prime of the largest primorial, showing thereby that the assumption is false that there is a largest prime number, N, that is the last twin prime.
- In Part 6, we present the proof of the Goldbach conjecture, in which we position an even number between two primorials, and then show that a Goldbach solution must exist for the even number.

## Part 1: Foundation

The Sieve of Eratosthenes
Eratosthenes (276-194 B.C.), born approximately 50 years after Euclid, was a contemporary of Archimedes. He created a sieve method to generate all prime numbers less than or equal to a given integer value N. Starting with 2, all integer multiples of 2 are struck from the list up to N; then the remaining multiples of 3 are removed, then the remaining multiples of 5 (4 being skipped since it is a multiple of 2), etc. The result is that once the multiples of the largest prime less than the square root of N have been eliminated in the sieve, there are no more values to be struck from the list, and what remains are all of the primes less than or equal to N. [5]

> Definition: Let N be a positive integer greater than or equal to 4. Then the set of seed primes in N equals the set of primes less than or equal to the square root of N.

For example, if N = 210 and the $\sqrt{N}$ = 14.49, then the set of seed primes = {2, 3, 5, 7, 11, 13}, where 13 is the largest prime less than or equal to 14.49.

The following are observations about the sieve of Eratosthenes:
- If N is the largest integer to be used to determine the elements in the sieve, then the largest prime used to generate the sieve is ≤ $\sqrt{N}$. Let $P_N$ equal the largest prime less than or equal to $\sqrt{N}$.
- As N increases to very large numbers, the ratio N / $\sqrt{N}$ increases accordingly (though more slowly). Similarly, for $P_N$ equal to the largest prime ≤ $\sqrt{N}$ the ratio N / $P_N$ increases to very large numbers as N increases to very large numbers, as well as for the ratio of N to each prime less than $P_N$.
- For primes $P_r$ ≤ $P_N$, as the ratio of N / $P_r$ increases, it permits more multiples of $P_r$ to be generated in the sieve. The result is that the increase in the number of primes slows as N increases while the corresponding increase in the number of composites grows.
- As each prime $P_r$ ≤ $P_N$ generates new composites, its first new composite occurs at $P_r^2$ and each subsequent new composite of $P_r$ is a multiple of $P_r$ and one or more primes greater than or equal to $P_r$. Consequently, each $P_r$ produces fewer new multiples < N than its predecessor in the sieve.

Congruence Classes

If S is a divisor of X – r, then we say that X is congruent to r modulus S, and we write

$$X \equiv r \text{ Mod } S$$

and refer to r as the residue (or remainder) of X Mod S. If 0 ≤ r ≤ S, then r is the least residue of X Mod S. Further, two numbers, X and Y, are congruent if they have the same residue under Mod S. That is, if X ≡ r Mod S and Y ≡ r Mod S, then we say that X and Y are congruent under Mod S. [6] For example, since 47 ≡ 2 Mod 5 and 7 ≡ 2 Mod 5, we say that 47 and 7 are congruent under Mod 5 because they have the same modular residue when divided by 5.

A congruence class Mod S is the class of all numbers congruent to a given residue Mod S, and every member of the class, denoted by [r], is called a *representative* of the class. Further, for each S, there are S classes, represented by the integers 0, 1, 2, …S-1. For example, for S = 5, there are five congruence classes: [0], [1], [2], [3] and [4]. In our example, 47 and 7 are members of congruence class [2] for Mod 5.

- For each modulo P where P is a prime number, congruence class [0] has only one prime number representative: the prime number P itself. All other members of [0] are multiples of P.
- All integers in class [1] are of the form n*P + 1, for integers n ≥ 0. Similarly, all integers in [2] are of the form n*P + 2, for integers n ≥ 0, and so forth for successively higher integers.
- In Mod 3, the prime numbers > 3 are in either congruence class [1] or [2]. Similarly, under Mod 5, the prime numbers not equal to 5 are distributed over the four non-zero congruence classes. In fact, as P gets bigger, the primes ≠ P are distributed over P-1 congruence classes (all those ≠ [0]), with the same separation (or gap) as their original differences, mod P (i.e. $P_r - P_j \equiv [r_i] - [r_j]$ Mod P).

Modular Signatures of the Positive Integers

Table 1 below shows the seed primes (2, 3 and 5) for the even number 30 = 2*3*5 which the Sieve of Eratosthenes uses to generate the prime values larger than 5 and less than 30. The table also shows the "signature" of each integer, consisting of its residues under modular arithmetic for each seed prime, ordered from the smallest seed prime to the largest. Because of the Fundamental Theorem of Arithmetic [7], each positive integer has a unique "signature" of modular residues under all of the prime numbers.

> Definition: The <u>modular signature</u> of a positive integer, Z, is its unique ordered (on prime numbers $P_r$ in ascending order starting with 2) sequence of modular residues, $(s_1, s_2, …, s_N)$, under the prime numbers less than or equal to some maximum seed prime, $P_N$, where
>
> $$Z \equiv s_i \text{ Mod } P_r \text{ for each } 2 \leq P_r \leq P_N$$

For example, the modular signature of the positive integer (and prime number) 13 under seed primes up to 19 is: (13 Mod 2, 13 Mod 3, 13 Mod 5, 13 Mod 7, 13 Mod 11, 13 Mod 13, 13 Mod 17, 13 Mod 19) which is equal to this modular signature: (1, 1, 3, 6, 2, 0, 13, 13).

Of particular interest is the fact that the modular residue $s_i$ of any integer, Z, equals Z for all seed primes greater than Z. That ensures that the modular signatures of every positive integer will be unique when seed primes are used that are greater than Z. For example, the modular signature of the integer 7 is the sequence (1, 1, 2, 0, 7, 7, 7, 7, 7, 7) under the first ten seed primes 2, 3, 5, 7, 11, 13, 17, 19, 23 and 29.

Table 1 below shows all possible combinations of modular residues for E = 30, using the factors of 30 (2, 3 and 5) as seed primes, generating a complete cycle of all of the combinations from 1 to 30 of the seed primes' modular residues.

Starting at 31, each of the seed primes 2, 3 and 5 will repeat its residue cycle <u>in the same order</u> for every increment of 30, so obviously more seed primes > 5 are needed to maintain unique modular signatures for integers > 30. Consequently, modular signature for integers > 30 will have to expand to include the modular residue(s) of one or more additional seed primes so as to eliminate the possibility of duplicate signatures.

| Integer | Prime | Modular Signatures Seed Primes | | | Twin Prime |
|---|---|---|---|---|---|
| | | 2 | 3 | 5 | |
| 1 | | 1 | 1 | 1 | |
| 2 | X | 0 | 2 | 2 | |
| 3 | X | 1 | 0 | 3 | |
| 4 | | 0 | 1 | 4 | |
| 5 | X | 1 | 2 | 0 | X |
| 6 | | 0 | 0 | 1 | |
| 7 | X | 1 | 1 | 2 | X |
| 8 | | 0 | 2 | 3 | |
| 9 | | 1 | 0 | 4 | |
| 10 | | 0 | 1 | 0 | |
| 11 | X | 1 | 2 | 1 | |
| 12 | | 0 | 0 | 2 | |
| 13 | X | 1 | 1 | 3 | X |
| 14 | | 0 | 2 | 4 | |
| 15 | | 1 | 0 | 0 | |
| 16 | | 0 | 1 | 1 | |
| 17 | X | 1 | 2 | 2 | |
| 18 | | 0 | 0 | 3 | |
| 19 | X | 1 | 1 | 4 | X |
| 20 | | 0 | 2 | 0 | |
| 21 | | 1 | 0 | 1 | |
| 22 | | 0 | 1 | 2 | |
| 23 | X | 1 | 2 | 3 | |
| 24 | | 0 | 0 | 4 | |
| 25 | | 1 | 1 | 0 | |
| 26 | | 0 | 2 | 1 | |
| 27 | | 1 | 0 | 2 | |
| 28 | | 0 | 1 | 3 | |
| 29 | X | 1 | 2 | 4 | |
| 30 | | 0 | 0 | 0 | |

Note: Only the larger of the pair is listed for twin primes

Table 1

In the next section, we will define the maximum seed prime for all modular signatures used in this paper, with the requirement that the maximum seed prime be large enough to eliminate the possibility of duplicate modular signatures.

Primorials and Modular Signatures
The product of consecutive prime numbers starting with 2 and up to a largest prime, P, is called a primorial and is denoted by P#. [8]  In this paper we will use primorials to frame modular signatures to ensure that the modular signatures contain enough residue elements to show that they are unique for each positive integer.

Note the following about Table 1 and primorial 30 = 2*3*5:
1. 30 is the only primorial for which its largest seed prime (5) also is the largest prime factor of 30.
2. Since there are two modular residues for seed prime 2, three residues for seed prime 3 and five residues for seed prime 5, there are 2*3*5 = 30 possible combinations of the residues to produce the modular signatures of the integers from 1 to 30.  Thus, the seed primes of 30 generate all possible combinations of modular signatures of the residues under seed primes 2, 3 and 5 in primorial 30.
3. The three seed primes of 30 each have modular signatures containing one 0, which is the residue for each respective seed prime when divided by itself.  Otherwise, a prime number greater than 5 and less than 30 has no zeros in its modular signature.
4. While there will be other integers greater than a seed prime that may have just one 0 in their modular signature (e.g. the number 25), they are composites because they are the product of more than a single power of the seed prime and the number 1 (e.g. $25 = 5^2 * 1$).

5. We can use the following formula to count the number of primes in primorial 30: 3 + [(2-1)*(3-1)*(5-1)] -1 = 10, where 3 equals the number of seed primes, [(2-1)*(3-1)*(5-1)] equals the number of combinations of non-zero residues of the three seed primes and -1 is included to subtract the number 1 from the positive combinations because 1 is not a prime number by definition.

To proceed, then, we will need to have a consistent method for determining the maximum seed prime required to produce the unique modular signature for a positive integer.

> Definition: Let N be a positive integer. Then the <u>maximum seed prime</u> for N's modular signature will be the largest prime less than or equal to the square root of the smallest primorial, E, greater than or equal to N.

For example, if N = 68, then E = 7# = 210 = 2*3*5*7 is the smallest primorial greater than or equal to 68. Since 13 is the largest prime less than $\sqrt{210}$, the modular signature for 68 will use seed primes 2, 3, 5, 7, 11 and 13, and every positive integer less than primorial 210 (including 68) will have a unique modular signature under the seed primes up to and including 13.

Whereas the Sieve of Eratosthenes uses the $\sqrt{N}$ to identify the maximum prime that will be used to generate primes less than or equal to a positive integer N, we use an extension of that concept with primorials to determine the seed primes that will be used to produce the unique modular signature of each positive integer.

**Theorem 1**. Every positive integer, N, has a unique modular signature under the prime numbers less than or equal to the maximum seed prime for N's modular signature, $P_Z$.

Proof.
1. Assume that M and N are two different positive integers under P#, the smallest primorial greater than both, where M ≠ N.
2. Let Q = the set of seed primes of primorial P#, with the largest element of Q equal to $P_Z$.
3. Let X = ($r_1$, $r_2$, $r_3$, …$r_Z$) be the modular signature of M under P#.
4. Let Y = ($s_1$, $s_2$, $s_3$,…$s_Z$) be the modular signature of N under P#.
5. Assume X = Y. That is, for every element of X and corresponding element of Y, $r_i \equiv s_i$ (mod $p_i$) for each $p_i \in Q$.
6. By the <u>Chinese Remainder Theorem</u> [9], the set of equations U ≡ $r_i$ ≡ $s_i$ (mod $p_i$) has a unique solution U modulo P#. Since U is a unique solution, U = M and U = N, and M = N.

**Theorem 2**. Let N be a positive integer and let A be the smallest primorial greater than or equal to N, where N is not a seed prime of A. If N has a modular signature under the seed primes of A in which all of the modular residues are non-zero, then N will be a prime number.

Proof.
1. Let P = the set of seed primes of A that are used to create the modular signature of N.
2. Let $P_N$ be the largest element in P; therefore, $P_N$ is the largest prime < $\sqrt{A}$. Since N ≤ A, $\sqrt{N} \leq \sqrt{A}$ and $P_N$ is greater than or equal to the largest prime ≤ $\sqrt{N}$.
3. Since the modular residue, $s_i$, for every $P_r \in P$ is non-zero in the modular signature of N, N is pairwise co-prime to every element of P.
4. Since N is pairwise co-prime to every seed prime of A, no prime number less than $\sqrt{N}$ is a factor of N. Then N is only divisible by 1 and N, so by definition, N is a prime number.

### Core Versus Non-Core Seed Primes

Because primorials larger than 30 will have seed primes that are greater than the largest prime factor of the primorial, it will be useful to distinguish between the seed primes which are factors of the primorial and those which are not factors of the primorial.

> Definition: Let N be a primorial greater than 30, and let $P_N$ equal the largest prime factor of N and $P_{MAX}$ equal the largest prime ≤ $\sqrt{N}$. Then the core seed primes for N will be the primes less than or equal to the largest prime factor of N, $P_N$. The non-core seed primes of N are the remaining seed primes of N, those which are greater than $P_N$ and less than or equal to $P_{MAX}$.

For example, for primorial 7# = 210 = 2*3*5*7, the core seed primes are the set {2, 3, 5, 7} and the non-core seed primes are the set {11, 13}. By definition, the two sets are disjoint.

As we saw with primorial 30, the core seed primes of a primorial will generate all possible combinations of their residue values in the modular signatures between 1 and the primorial, inclusive. Further, using the combinatorial Multiplication Principle [10] in a primorial, there will be X combinations of non-zero residues of the N core seed primes, where $X = \prod_{r=1}^{N}(Pr - 1)$ for all primes $P_r \leq P_N$, the largest prime factor of the primorial (the largest core seed prime). All of those combinations will contain residues greater than 0 in their modular signatures for the core seed primes, and all will apply to odd integers.

> Definition: Let M be a positive odd integer less than a primorial, E. Then M is a potential prime under E if, in its modular signature under E, M has a non-zero modular residue for every core seed prime of E.

Note that the formula for the number of potential primes in a primorial, $X = \prod_{r=1}^{N}(Pr - 1)$ for all primes $P_r \leq P_N$, is equivalent to Euler's Totient Function [11] when the Totient Function is applied to the primorial. However, some of the potential primes in the primorial will become composites, because the modular signatures of some of the potential primes will contain residues equal to zero for one or more of the non-core seed primes. For example, as shown in Table 2 below, the modular signature of integer 2291 has all positive residues in its five core seed primes under primorial 2310, but has a zero residue for non-core seed prime 29, because 2291 = 29*79. Therefore, 2291 is a potential prime number under primorial 2310, but is converted to a composite by a zero residue in its modular signature generated by a non-core seed prime (those greater than $P_r$ = 11 and less than or equal to 47).

|         |    |    |    |    |    | Mod Residues For Seed Primes |    |    |    |    |    |    |    |    |    |
|---------|----|----|----|----|----|----|----|----|----|----|----|----|----|----|----|
|         | Core Seed Primes | | | | | Non-Core Seed Primes | | | | | | | | | |
| Integer | 2  | 3  | 5  | 7  | 11 | 13 | 17 | 19 | 23 | 29 | 31 | 37 | 41 | 43 | 47 |
|         |    |    |    |    |    |    |    |    |    |    |    |    |    |    |    |
| 2,291   | 1  | 2  | 1  | 2  | 3  | 3  | 13 | 11 | 14 | 0  | 28 | 34 | 36 | 12 | 35 | < Signature

Table 2

Further, if none of the potential prime numbers had a zero residue with a non-core seed prime, then the number of prime numbers in primorial 2310 would = 5 + ([2-1]*[3-1]*[5-1]*[7-1]*[11-1]) -1 = 5+480-1 = 484 instead of 343, the actual number of primes less than 2310. The reduction of 141 (29.1%) was due to the impact of the zero residues generated by the 10 non-core seed primes (13 to 47) when they were included in the modular signatures of the potential prime numbers.

### Primorial Stacking

Each primorial > 6 is the product of its predecessor primorial and the next larger prime number in the sequence of primes. Thus, the modular residues under core seed primes 2, 3 and 5 in primorial 5# = 30 = 2*3*5 are repeated 7 times in the modular signatures of primorial 7# = 210 = 30*7 and 77=7*11 times in the modular signatures under primorial 11# = 2310. Similarly, the modular residues under core seed primes 2, 3, 5 and 7 in primorial 210 are repeated 11 times in the modular signatures of primorial 11# = 2310 = 210*11. Of course, each of the modular signatures becomes unique when we include their non-core seed prime residues in each integer's modular signature under primorial 2310's seed primes.

Table 3 below illustrates primorial stacking of 5# = 30 = 2*3*5 within primorial 7# = 210 = 30*7. Under primorial 210, the largest prime less than $\sqrt{210}$ is 13, so the core seed primes of 210 (2, 3, 5 and 7) are joined by two non-core seed primes, 11 and 13.

To illustrate the repetition of the modular signatures of primorial 30 in primorial 210, the odd integers 1 through 195 are grouped in Table 3 in congruence classes Mod 30, each with 7 elements. That is, the class [1] Mod 30 = {1, 31, 61, 91, 121, 151, 181} for odd integers less than 210. Similar groupings are shown for congruence classes [3], [5], [7], [9], [11], [13] and [15] Mod 30 over the set of odd integers less than or equal to 210, and could be shown for the remaining seven congruence classes for Mod 30 (17, 19, 21, 23, 25, 27 and 29) with similar results covering the remaining odd integers through 209.

| | | Modular Signatures of Odd Integers (1-195) For Primorial 210, Grouped By Mod 30 Congruence Classes | | | | | | | | | | | | |
|---|---|---|---|---|---|---|---|---|---|---|---|---|---|---|
| | | Seed Primes | | | | | | | | Seed Primes | | | | |
| | | Core Seed Primes | | | | Non-Core | | | | Core Seed Primes | | | | Non-Core | |
| Odd | Primes | 2 | 3 | 5 | 7 | 11 | 13 | Odd | Primes | 2 | 3 | 5 | 7 | 11 | 13 |
| 1 | | 1 | 1 | 1 | 1 | 1 | 1 | 9 | | 1 | 0 | 4 | 2 | 9 | 9 |
| 31 | X | 1 | 1 | 1 | 3 | 9 | 5 | 39 | | 1 | 0 | 4 | 4 | 6 | 0 |
| 61 | X | 1 | 1 | 1 | 5 | 6 | 9 | 69 | | 1 | 0 | 4 | 6 | 3 | 4 |
| 91 | | 1 | 1 | 1 | 0 | 3 | 0 | 99 | | 1 | 0 | 4 | 1 | 0 | 8 |
| 121 | | 1 | 1 | 1 | 2 | 0 | 4 | 129 | | 1 | 0 | 4 | 3 | 8 | 12 |
| 151 | X | 1 | 1 | 1 | 4 | 8 | 8 | 159 | | 1 | 0 | 4 | 5 | 5 | 3 |
| 181 | X | 1 | 1 | 1 | 6 | 5 | 12 | 189 | | 1 | 0 | 4 | 0 | 2 | 7 |
| | | | | | | | | | | | | | | | |
| 3 | X | 1 | 0 | 3 | 3 | 3 | 3 | 11 | X | 1 | 2 | 1 | 4 | 0 | 11 |
| 33 | | 1 | 0 | 3 | 5 | 0 | 7 | 41 | X | 1 | 2 | 1 | 6 | 8 | 2 |
| 63 | | 1 | 0 | 3 | 0 | 8 | 11 | 71 | X | 1 | 2 | 1 | 1 | 5 | 6 |
| 93 | | 1 | 0 | 3 | 2 | 5 | 2 | 101 | X | 1 | 2 | 1 | 3 | 2 | 10 |
| 123 | | 1 | 0 | 3 | 4 | 2 | 6 | 131 | X | 1 | 2 | 1 | 5 | 10 | 1 |
| 153 | | 1 | 0 | 3 | 6 | 10 | 10 | 161 | | 1 | 2 | 1 | 0 | 7 | 5 |
| 183 | | 1 | 0 | 3 | 1 | 7 | 1 | 191 | X | 1 | 2 | 1 | 2 | 4 | 9 |
| | | | | | | | | | | | | | | | |
| 5 | X | 1 | 2 | 0 | 5 | 5 | 5 | 13 | X | 1 | 1 | 3 | 6 | 2 | 0 |
| 35 | | 1 | 2 | 0 | 0 | 2 | 9 | 43 | X | 1 | 1 | 3 | 1 | 10 | 4 |
| 65 | | 1 | 2 | 0 | 2 | 10 | 0 | 73 | X | 1 | 1 | 3 | 3 | 7 | 8 |
| 95 | | 1 | 2 | 0 | 4 | 7 | 4 | 103 | X | 1 | 1 | 3 | 5 | 4 | 12 |
| 125 | | 1 | 2 | 0 | 6 | 4 | 8 | 133 | | 1 | 1 | 3 | 0 | 1 | 3 |
| 155 | | 1 | 2 | 0 | 1 | 1 | 12 | 163 | X | 1 | 1 | 3 | 2 | 9 | 7 |
| 185 | | 1 | 2 | 0 | 3 | 9 | 3 | 193 | X | 1 | 1 | 3 | 4 | 6 | 11 |
| | | | | | | | | | | | | | | | |
| 7 | X | 1 | 1 | 2 | 0 | 7 | 7 | 15 | | 1 | 0 | 0 | 1 | 4 | 2 |
| 37 | X | 1 | 1 | 2 | 2 | 4 | 11 | 45 | | 1 | 0 | 0 | 3 | 1 | 6 |
| 67 | X | 1 | 1 | 2 | 4 | 1 | 2 | 75 | | 1 | 0 | 0 | 5 | 9 | 10 |
| 97 | X | 1 | 1 | 2 | 6 | 9 | 6 | 105 | | 1 | 0 | 0 | 0 | 6 | 1 |
| 127 | X | 1 | 1 | 2 | 1 | 6 | 10 | 135 | | 1 | 0 | 0 | 2 | 3 | 5 |
| 157 | X | 1 | 1 | 2 | 3 | 3 | 1 | 165 | | 1 | 0 | 0 | 4 | 0 | 9 |
| 187 | | 1 | 1 | 2 | 5 | 0 | 5 | 195 | | 1 | 0 | 0 | 6 | 8 | 0 |

Table 3

In each congruence class, the modular signatures of each member of the class are identical under the core seed primes for 30 = 2*3*5, and differ under seed prime 7. Under primorial stacking, within primorial 210 we would expect to see 7 iterations of the modular signatures of each of the positive integers less than or equal to primorial 30, with the signatures differing in seed primes > 5. For example, within congruence class [1] = {1, 31, 61, …181} the seven integers have identical modular signatures under core seed primes 2, 3 and 5, and the residues of core seed prime 7 are added to each member of [1] starting with 1 to help create a unique modular signature for each integer under primorial 210.

Further, the residues for non-core seed prime 11 increase by 8 = 30 Mod 11 between corresponding cycles of primorial 30 in each congruence class. That is, the residue for 31 Mod 11 equals 1 Mod 11 + 30 Mod 11 = 1 + 8 = 9. Similarly, the residues for non-core seed prime 13 increase by 4 = 30 Mod 13 between cycles of primorial 30 in each congruence class (e.g. between 1 and 31 Mod 13). Were the modular signatures expanded to include seed primes 17, 19, etc., with the next primorial, then the residues Mod $P_r$ of an integer in each congruence class would increase similarly by $s_i$ = 30 Mod $P_r$ with each cycle of 30.

The same can be said for the other congruence classes Mod 30: the members of each congruence class share identical modular signatures under primorial 30's core seed primes 2, 3 and 5, and differ under core seed prime 7, with one member in each congruence class receiving the 0 residue from core seed prime 7 in the middle of the cycle of the residues of 7 until 7's residue cycle starts again after 210.

That pattern will continue when primorial 210 is stacked in the next larger primorial, 11# = 2310 = 210*11. The modular signature of 210 will contain zeros for its core seed primes (2, 3, 5, and 7), causing 211 to restart the cycles of 210's core seed primes (2, 3, 5, and 7) at residue 1 in each of its multiples within primorial 2310. The core seed primes of 210 (2, 3, 5 and 7) are extended under primorial 2310 with an additional core seed prime (11) and ten non-core seed primes (13, 17, 19, 23, 29, 31, 37, 41, 43, 47). Since the modular signature of 210 is (0,0,0,0) under the core seed primes 2, 3, 5 and 7, when 210 is added to any integer less than 2310, the modular signature of the smaller number will be repeated for the core seed primes of 210. This is illustrated in Table 4 below with two congruence classes of 210 at opposite ends of the primorial: [1] and [209].

| | | | Mod Signatures Under Primorial 2310 = 210*11 | | | | | | | | | | | | | |
|---|---|---|---|---|---|---|---|---|---|---|---|---|---|---|---|---|
| | | | Core Seed Primes | | | | | Non-Core Seed Primes | | | | | | | | |
| | | | 2 | 3 | 5 | 7 | 11 | 13 | 17 | 19 | 23 | 29 | 31 | 37 | 41 | 43 | 47 | Prime |
| Count | Odds | Delta | | | | | | | | | | | | | | | | |
| 1 | 1 | | 1 | 1 | 1 | 1 | 1 | 1 | 1 | 1 | 1 | 1 | 1 | 1 | 1 | 1 | 1 | |
| 2 | 211 | 210 | 1 | 1 | 1 | 1 | 2 | 3 | 7 | 2 | 4 | 8 | 25 | 26 | 6 | 39 | 23 | X |
| 3 | 421 | 210 | 1 | 1 | 1 | 1 | 3 | 5 | 13 | 3 | 7 | 15 | 18 | 14 | 11 | 34 | 45 | X |
| 4 | 631 | 210 | 1 | 1 | 1 | 1 | 4 | 7 | 2 | 4 | 10 | 22 | 11 | 2 | 16 | 29 | 20 | X |
| 5 | 841 | 210 | 1 | 1 | 1 | 1 | 5 | 9 | 8 | 5 | 13 | 0 | 4 | 27 | 21 | 24 | 42 | |
| 6 | 1051 | 210 | 1 | 1 | 1 | 1 | 6 | 11 | 14 | 6 | 16 | 7 | 28 | 15 | 26 | 19 | 17 | X |
| 7 | 1261 | 210 | 1 | 1 | 1 | 1 | 7 | 0 | 3 | 7 | 19 | 14 | 21 | 3 | 31 | 14 | 39 | |
| 8 | 1471 | 210 | 1 | 1 | 1 | 1 | 8 | 2 | 9 | 8 | 22 | 21 | 14 | 28 | 36 | 9 | 14 | X |
| 9 | 1681 | 210 | 1 | 1 | 1 | 1 | 9 | 4 | 15 | 9 | 2 | 28 | 7 | 16 | 0 | 4 | 36 | |
| 10 | 1891 | 210 | 1 | 1 | 1 | 1 | 10 | 6 | 4 | 10 | 5 | 6 | 0 | 4 | 5 | 42 | 11 | |
| 11 | 2101 | 210 | 1 | 1 | 1 | 1 | 0 | 8 | 10 | 11 | 8 | 13 | 24 | 29 | 10 | 37 | 33 | |
| 1 | 209 | | 1 | 2 | 4 | 6 | 0 | 1 | 5 | 0 | 2 | 6 | 23 | 24 | 4 | 37 | 21 | |
| 2 | 419 | 210 | 1 | 2 | 4 | 6 | 1 | 3 | 11 | 1 | 5 | 13 | 16 | 12 | 9 | 32 | 43 | X |
| 3 | 629 | 210 | 1 | 2 | 4 | 6 | 2 | 5 | 0 | 2 | 8 | 20 | 9 | 0 | 14 | 27 | 18 | |
| 4 | 839 | 210 | 1 | 2 | 4 | 6 | 3 | 7 | 6 | 3 | 11 | 27 | 2 | 25 | 19 | 22 | 40 | X |
| 5 | 1049 | 210 | 1 | 2 | 4 | 6 | 4 | 9 | 12 | 4 | 14 | 5 | 26 | 13 | 24 | 17 | 15 | X |
| 6 | 1259 | 210 | 1 | 2 | 4 | 6 | 5 | 11 | 1 | 5 | 17 | 12 | 19 | 1 | 29 | 12 | 37 | X |
| 7 | 1469 | 210 | 1 | 2 | 4 | 6 | 6 | 0 | 7 | 6 | 20 | 19 | 12 | 26 | 34 | 7 | 12 | |
| 8 | 1679 | 210 | 1 | 2 | 4 | 6 | 7 | 2 | 13 | 7 | 0 | 26 | 5 | 14 | 39 | 2 | 34 | |
| 9 | 1889 | 210 | 1 | 2 | 4 | 6 | 8 | 4 | 2 | 8 | 3 | 4 | 29 | 2 | 3 | 40 | 9 | X |
| 10 | 2099 | 210 | 1 | 2 | 4 | 6 | 9 | 6 | 8 | 9 | 6 | 11 | 22 | 27 | 8 | 35 | 31 | X |
| 11 | 2309 | 210 | 1 | 2 | 4 | 6 | 10 | 8 | 14 | 10 | 9 | 18 | 15 | 15 | 13 | 30 | 6 | X |
| | 30 | | 0 | 0 | 0 | 2 | 8 | 4 | 13 | 11 | 7 | 1 | 30 | 30 | 30 | 30 | 30 | |
| | 210 | | 0 | 0 | 0 | 0 | 1 | 2 | 6 | 1 | 3 | 7 | 24 | 25 | 5 | 38 | 22 | |
| | 2310 | | 0 | 0 | 0 | 0 | 0 | 9 | 15 | 11 | 10 | 19 | 16 | 16 | 14 | 31 | 7 | |

Table 4

Just as the members of the 15 odd congruence classes under Mod 30 had identical modular residues under the core seed primes of 30, when 30 is stacked in 210 and the odd integers are grouped into 105 congruence classes under Mod 210, all of the modular residues of the members of a congruence class are identical under the core seed primes of both 30 and 210, and then differ in the remaining seed primes that make up their modular signatures. As the primorials are stacked further, this property is repeated with the addition of each core seed prime. In short, there is a pattern to the modular signatures in each congruence class of a primorial that is repeated in every larger primorial, no matter how large.

This pattern of dispersion of the modular signatures within successive primorials, and the dispersion of the integers with non-zero residues under the core seed primes of a primorial, will be relevant later when we discuss the dispersion of potential primes, potential twin primes and potential Goldbach solutions within primorials.

New Composites of Non-Core Seed Primes
Returning to the Sieve of Eratosthenes, we saw that the first <u>new</u> composite (i.e. those not already generated by a smaller prime) for a seed prime $P_r$ occurred at $P_r^2$. We also noted that each new composite thereafter for $P_r$ is generated as the product of seed prime $P_r$ and a combination of primes greater than or equal to $P_r$. Within primorial 2310, for example, the non-core seed prime 13's first new composite is 169 = 13*13. Its next new composite is 221 = 13*17. 13's largest new composite in primorial 2310 is 2249 = 13*173. There are larger composites in primorial 2310 that have 13 as a factor (e.g. 2275 = $5^2$*7*13), but those composites were generated first by primes smaller than 13. In fact, when we partition the positive integers into primorials, and examine only the composites generated by the primorial's non-core seed primes, we can see that all of the <u>prime factors</u> of the "new composite" integers generated by the non-core seed primes in a primorial are contained within the next <u>smaller</u> primorial.

**Theorem 3**: Under modular signatures for a primorial E, the first new composite (i.e. not already generated by a smaller prime) generated by a seed prime, $P_r$, of E, will occur at $P_r^2$. All composites less than $P_r^2$ in E will contain a prime factor less than $P_r$.

Proof.
1. Let A be any composite positive integer less than $P_r^2$ in E.
2. By definition, A has a prime divisor X with 1 < X < A, thus A = XY.
3. Either $X \leq \sqrt{A}$ or $Y \leq \sqrt{A}$. Otherwise, if $X > \sqrt{A}$ and $Y > \sqrt{A}$, then $XY > \sqrt{A} * \sqrt{A} = A$.
4. Therefore, if A is a composite number less than $P_r^2$, it has a prime divisor less than $P_r$, since $\sqrt{A} < P_r$.

**Corollary 3.1**: Let M be a primorial and let N be the smallest primorial larger than M. If $C_i$ is a new composite of N generated by a non-core seed prime, $P_r$, of N, then the other prime factors of $C_i$ will be < M.

Proof.
1. Let $P_N$ be the largest prime factor of N such that $N = M * P_N$.
2. Since $P_r$ is a non-core seed prime of N, $P_N < P_r$.
3. By <u>Theorem 3</u>, the first composite in N of a non-core seed prime, $P_r$, will be $P_r^2$.
4. Then $P_r^2 \leq C_i < N$.
5. By definition, $C_i$ has a divisor X with $P_r^2 \leq P_r * X < N$.
6. Then $P_r \leq X < (N / P_r)$ and $P_r \leq X < ((M*P_N) / P_r)$.
7. Since $P_N < P_r$, X < M, and any prime factor of X will be less than M.

**Corollary 3.2**: Let N be a potential prime number under the core seed primes of primorial E, with N not equal to any non-core seed prime of E. Let Q equal the set of non-core seed primes of E. If $N < P_N^2$, where $P_N$ is the smallest element of Q, then the modular residue of N Mod $P_S$ will be greater than zero for all $P_S \in Q$.

Proof:
1. By definition of a potential prime, the modular signature of N will not contain a zero residue for any of the core seed primes of E.
2. Assume N is a composite of $P_N$, with N < $P_N^2$.
3. By <u>Theorem 3</u>, since N is a composite < $P_N^2$, it must contain a prime factor less than $P_N$. Since $P_N$ is the smallest non-core seed prime of E, the prime factor less than $P_N$ must be a core seed prime of E.
4. However, since N is a potential prime of E, it cannot have a prime factor equal to a core seed prime of E.
5. Then N ≠ 0 Mod $P_N$ and the modular residue of N Mod $P_N$ is greater than zero.
6. Let $P_S \in Q$ be the first element of Q which is a factor of N.
7. Since N < $P_N^2$, and $P_N$ is the smallest element of Q, N = X * $P_S$ for some X < $P_N$. However, since N is a potential prime of E and $P_N$ is not a factor of N, there cannot exist a prime number X < $P_N$ that is a factor of N. Therefore, N Mod $P_S$ > 0 for all $P_S \in Q$.

<u>Formula for Computing the Number of Primes in a Primorial</u>
Generally, the Sieve of Eratosthenes takes this form for the number of primes less than a positive integer, N:

$$\pi(N) = N - 1 + n(SP) - n(\text{Composites of SP}) \quad (1)$$

where 1 is subtracted since 1 is not a prime by definition, SP equals the set of seed primes for N, n(SP) is the number of seed primes and n(Composites of SP) is the number of unique elements in the set of composites less than or equal to N that are generated by the seed primes. For example, $\pi(100) = 100 - 1 + 4 - (117 - 45 + 6 - 0) = 25$.

Further, since the seed primes can be divided into two discrete sets of core and non-core seed primes, the composites of the seed primes can be divided into two sets: A = the composites generated by the core seed primes and B = the new composites generated by the non-core seed primes. By the combinatorial Inclusion-Exclusion Principle [12],

$$n(\text{Composites of SP}) = n(A) + n(B) - n(A \cap B) \quad (2)$$

Define P as the set of prime numbers and $P_c$ as the set of core seed primes for N, and $P_{nc}$ as the set of non-core seed primes for N, then create the set of new composites in A by multiplying each element $p_r$ of $P_c$ by primes greater than or equal to $p_r$, such that their product is less than or equal to N. Similarly, create the set of new composites in B by multiplying each element $p_j$ of $P_{nc}$ only by primes greater than or equal to $p_j$, such that their product is less than or equal to N. Then the number of elements in (A ∩ B), n(A ∩ B), equals 0 since the composites in B will have no factors from among the core seed primes and therefore will not equal any of the composites in A (i.e. A ∩ B is the empty set, and A and B are disjoint sets).

Now define an alternate form of the sieve of Eratosthenes <u>for primorials</u>. For primorial E,

$$\pi(E) = M + \left( \prod_{z=1}^{M}(Pz - 1) \right) - 1 - n(B) \quad (3)$$

where CSP is the set of core seed primes of E, M is the number of elements in CSP, Pz are the elements of CSP (non-core seed primes are counted within $\prod(Pz - 1)$), and B is the set of new composites generated by the non-core seed primes of E. The product of the core seed primes of E, $\prod(Pz - 1)$, will generate all possible combinations of the non-zero residues of the <u>core</u> seed primes in the modular signatures of the integers from 1 to E. Those are all potential prime numbers, including 1 (which is subtracted in the formula). However, some of those potential primes will become composites when the non-core seed primes add one or more [0] residue elements to their modular signatures. That was illustrated for 2291 in Table 2, because of the zero-residue added by 29 in 2291 = 29 x 79.

Therefore, we need to determine if there is a pattern to when the non-core seed primes eliminate potential prime numbers from the subset of potential prime numbers that could be twin primes or Goldbach solutions, knowing that each non-core seed prime generates new composites using prime factors greater than or equal to itself and less than the next smaller primorial. Establishing that pattern enables the proofs.

Modular Residue Cycles
Consider the behavior of the residues of the integers under modular arithmetic, extending the picture illustrated in Table 1 to include an additional seed prime, 7. As we move down the column of integers, incrementing each by adding 1 to the preceding value, we move through the modular remainders in a synchronized way.

Note that each seed prime operates a repeating cycle consisting of its residues through successive positive integers, all starting where N = 1. As we increment N by 1 (i.e. cycle through the integers), each mod residue increases by 1. And each seed prime, $P_r$ in turn, will complete a full cycle (equal in length to twice the prime for the odd integers) at a different point, because each possesses $P_r$ residues, with the difference in the ending integer between two primes determined by the size of the gap between the two primes.

It should now be obvious that there is a logic to when the prime numbers appear within the set of positive integers. A prime appears within a primorial when each residue of its modular signature under each smaller seed prime is non-zero, and its own modular residue equals 0 for the first time. That happens in the modular signature of a prime number when the repeating sequence of modular residues of the seed primes all fall on a non-zero value for each seed prime except for the prime itself. As noted earlier, the powers of a prime number (e.g. $7^2$ = 49) can also have a single 0 residue in their modular signature, but they will not be a prime because they will not be the first occurrence of that 0 residue among the integers for that seed prime (e.g. 7 is the first occurrence of the 0-residue mod 7, while 49 is not, as 49 is divisible by both 7 and 49).

Potential Prime and Non-Core Seed Prime Composite Patterns
In finding the number of primes in a primorial, the two main variable components are (1) the number of potential prime numbers generated by all possible combinations of the non-zero residues of the core seed primes (thereby eliminating the composites of the core seed primes), and (2) the number of new composites generated with a zero residue by each of the non-core seed primes in the modular signatures of each odd number, thereby converting some of the potential primes into composites.

Figure 1 below depicts graphically a view of the non-core seed prime composite data for primorial 30,030 along with the potential primes, by cycle of the largest seed prime (and largest non-core seed prime), 173. The core seed primes (2, 3, 5, 7, 11 and 13) generated a total of 5,760 potential primes among the odd integers, and the non-core seed primes (17 through 173) generated 2,517 unique new composites. Thus, the number of primes for primorial 30,030 = $\pi$ (30,030) = 6 + (5,760 – 1) – 2,517 = 3,248, where 6 is the number of core seed primes for primorial 30,030 and 1 is subtracted from 5,760 because it is included in the potential primes but is not a prime number by definition.

As expected because of the cycling of the modular residues in the modular signatures of the core seed primes, the number of potential primes is fairly constant across each cycle of the largest non-core seed prime (cycle length = 2 * 173 = 346), varying primarily because of the different gap sizes between primes causing an overlap with the cycle cutoff.

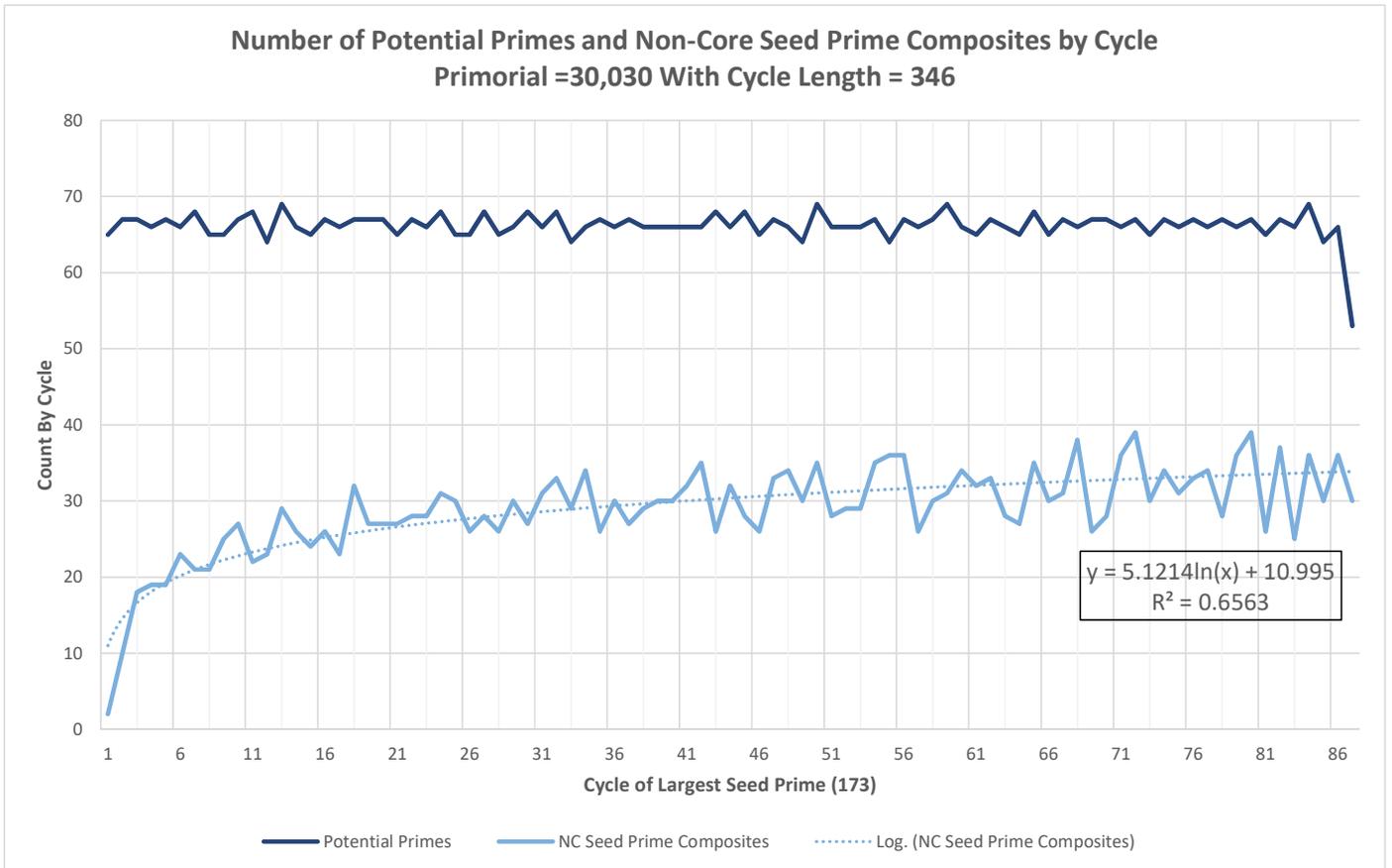

Figure 1

Therefore, to produce the gradually smaller number of primes per cycle characteristic of prime numbers, we expect the <u>cumulative</u> number of new non-core seed prime composites to increase gradually as the cycles unfold, starting at $P_r^2$ for each non-core seed prime $P_r$. That is exactly what happens as shown in Figure 1.

The consistent distribution of the potential primes across the span of the primorial is understandable in light of two factors: (1) primorial stacking and (2) the multiplicative property of the Euler Totient Function, $\varphi(n)$. [13]

**Theorem 4**: Let A equal a primorial whose largest prime factor (and therefore, its largest core seed prime) is $P_M$, and let B equal the next larger primorial with largest prime factor $P_N$ (i.e. B = A * $P_N$). Then the number of potential primes in B = $\varphi(B) = \varphi(A)\varphi(P_N)$, where $\varphi(P_N) = P_N - 1$. On average, $\varphi(B) / P_N$ potential primes of A will be repeated in each cycle of A within B.

Proof.
1. Since B is the next larger primorial after A, B = A*$P_N$ and the number of potential primes in B = $\varphi(B) = \varphi(A)\varphi(P_N)$ by the multiplicative property of the Euler Totient function.
2. Since $P_N$ is a prime number, $\varphi(P_N) = P_N - 1$.
3. Each potential prime in A contains a modular signature on the core prime factors in A, and each potential prime's modular signature in A is repeated $P_N$ times in the modular signatures of B. However, when the core seed primes of A are supplemented by the core seed prime $P_N$ of B, one of the modular signatures of a potential prime in A will contain a 0 residue from the congruence classes of $P_N$, leaving just $P_N-1$ non-zero residues for the potential primes of A in each cycle of B. Thus, on average, each cycle of A in B will have fewer potential primes in B than in A.

4. Then <u>on average</u> there will be φ(A) * [($P_N$ -1) / $P_N$] = φ(B) / $P_N$ potential primes of A repeated in each cycle of A within B.

**Part 2: Application of Modular Signatures to the Twin Primes Conjecture**

If the number of potential primes is fairly consistent across the cycle of the largest seed prime, can we expect the number of potential twin primes to be fairly consistent within primorial cycles? And can we expect the number of twin prime pairs to decline as the cycle count increases because of the increasing number of composites that will "knock out" potential twins from being true twins? And can we determine when the first potential twin prime will be converted to a "false" twin (or a composite) by a <u>non-core</u> seed prime's zero residue in the modular signature of the potential twin prime?

Table 5 below shows how the core seed primes generate potential twin primes among the potential primes for primorial 13# = 30030 = 2*3*5*7*11*13 = 2310*13. On average, the core seed primes generate 443 potential primes in each cycle of primorial 2310 in the primorial 30030. And on average, the seed primes generate a subset of 114 potential twins among those potential primes in each cycle of 2310 within 30030.

That result is expected given our earlier discussion of the multiplicative nature of the Euler Totient Function (Theorem 4), and as discussed below, the multiplicative nature of the function used to determine the number of potential twin primes. The reason: the repeating cycle of the modular residues of the seed primes in the modular signatures of each odd integer creates a consistent distribution of the non-zero residues in each cycle of smaller primorials within a larger primorial due to primorial stacking.

| | | | | Twin Prime Count By Cycle For Primorial 30030 | | | | | | |
|---|---|---|---|---|---|---|---|---|---|---|
| | | | | Cumulartive Number | | | | Number In Cycle | | |
| | Cycle - Primorial 2310 | | Potential | Potential | FALSE | TRUE | Potential | Potential | FALSE | TRUE |
| Count | End | Length | Primes | Twins | Twins | Twins | Primes | Twins | Twins | Twins |
| 1 | 2310 | | 443 | 113 | 47 | 66 | 443 | 113 | 47 | 66 |
| 2 | 4620 | 2310 | 887 | 228 | 112 | 116 | 444 | 115 | 65 | 50 |
| 3 | 6930 | 2310 | 1329 | 343 | 186 | 157 | 442 | 115 | 74 | 41 |
| 4 | 9240 | 2310 | 1772 | 455 | 266 | 189 | 443 | 112 | 80 | 32 |
| 5 | 11550 | 2310 | 2215 | 570 | 345 | 225 | 443 | 115 | 79 | 36 |
| 6 | 13860 | 2310 | 2658 | 684 | 430 | 254 | 443 | 114 | 85 | 29 |
| 7 | 16170 | 2310 | 3102 | 799 | 515 | 284 | 444 | 115 | 85 | 30 |
| 8 | 18480 | 2310 | 3545 | 913 | 593 | 320 | 443 | 114 | 78 | 36 |
| 9 | 20790 | 2310 | 3988 | 1028 | 677 | 351 | 443 | 115 | 84 | 31 |
| 10 | 23100 | 2310 | 4431 | 1141 | 753 | 388 | 443 | 113 | 76 | 37 |
| 11 | 25410 | 2310 | 4873 | 1255 | 846 | 409 | 442 | 114 | 93 | 21 |
| 12 | 27720 | 2310 | 5317 | 1370 | 932 | 438 | 444 | 115 | 86 | 29 |
| 13 | 30030 | 2310 | 5760 | 1484 | 1019 | 465 | 443 | 114 | 87 | 27 |
| | | | | | | | | | | |
| | Min | | | | | | 442 | 112 | 47 | 21 |
| | Max | | | | | | 444 | 115 | 93 | 66 |
| | Median | | | | | | 443 | 114 | 80 | 32 |
| | StdDev | | | | | | 0.62 | 0.95 | 11.29 | 11.10 |
| | | | | | | | | | | |
| Note: | "True Twins" in the first cycle excludes 5, 7, and 13; their predecessors had 0 in the modular residue (core seed primes) | | | | | | | | | |

Table 5

Further, consistent with Theorem 3, we would have expected the first new composite of primorial 13# = 30030 to convert a potential twin prime to a composite to appear when the first non-core seed prime, 17, equals $17^2$ = 289. However, 289 is not a potential twin prime because the odd numbers on either side of it (287 and 291) are not potential primes. Similarly, 323 = 17*19 is also not a potential twin prime for the same reason as 289. Consequently, the first new composite of a non-core seed prime to convert a potential twin prime to a composite in primorial 30030 is 361 = $19^2$. Therefore, because the set of potential twin primes is a subset of the set of potential primes, we can expect the first new composite of a non-core seed prime, $P_r$, to convert a potential twin prime to a composite no earlier than $P_r^2$, and it may happen later.

Predicting Potential Twin Primes From Modular Signatures

It is notable that there is a pattern in the modular residues that predicts when a potential prime can be a potential twin prime for two consecutive odd integers, O1 and O2 (with O1 < O2), where $s_i$ is the modular residue for a core seed prime, $P_r$:

- If $s_i = 0$ in the modular signature of O1 for a core seed prime, then O1 will not be a potential prime.
- If $s_i = P_r - 2$ in the modular signature of O1 for a core seed prime $P_r$, then O2 = O1 + 2 will not be a potential prime (and therefore, not a potential twin) since the second odd integer in the pair will be a composite.

The pattern is illustrated in Table 6 below which shows the modular signatures for ten consecutive odd integers within primorial 11# = 2310, where the largest core seed prime is 11 and the largest non-core seed prime (the largest prime less than $\sqrt{2310}$) is 47. Four of the odd integers are prime numbers, five are potential primes and just two are potential twin primes.

2237 and 2239 are both prime numbers because there are no residues equal to zero for either the core seed primes or the non-core seed primes in their modular signatures, and 2239 is a potential twin because both 2237 and 2239 are potential primes. 2241 is not a potential prime because it has a zero modular residue under Mod 3, not surprising because its predecessor, 2239, had 1 for the modular residue under Mod 3, and 1+2 = 0 Mod 3. For the same reason, 2247 and 2253 are not potential primes. Similarly, 2245 is not a potential prime because 2243 had 3 for the modular residue under Mod 5, and 3+2 = 0 Mod 5. For the same reason, 2255 is not a potential prime. The same type of residue combinations under Mod 7 caused 2245 to not be a potential prime, or a potential twin prime.

| Primes | Odds | Delta vs Prior | Core Seed Primes | | | | | Non-Core Seed Primes | | | | | | | | | Potential Prime | Potential Twin Prime |
|---|---|---|---|---|---|---|---|---|---|---|---|---|---|---|---|---|---|---|
| | | | 2 | 3 | 5 | 7 | 11 | 13 | 17 | 19 | 23 | 29 | 31 | 37 | 41 | 43 | 47 | | |
| X | 2237 | 2 | 1 | 2 | 2 | 4 | 4 | 1 | 10 | 14 | 6 | 4 | 5 | 17 | 23 | 1 | 28 | X | |
| X | 2239 | 2 | 1 | 1 | 4 | 6 | 6 | 3 | 12 | 16 | 8 | 6 | 7 | 19 | 25 | 3 | 30 | X | X |
| | 2241 | 2 | 1 | 0 | 1 | 1 | 8 | 5 | 14 | 18 | 10 | 8 | 9 | 21 | 27 | 5 | 32 | | |
| X | 2243 | 2 | 1 | 2 | 3 | 3 | 10 | 7 | 16 | 1 | 12 | 10 | 11 | 23 | 29 | 7 | 34 | X | |
| | 2245 | 2 | 1 | 1 | 0 | 5 | 1 | 9 | 1 | 3 | 14 | 12 | 13 | 25 | 31 | 9 | 36 | | |
| | 2247 | 2 | 1 | 0 | 2 | 0 | 3 | 11 | 3 | 5 | 16 | 14 | 15 | 27 | 33 | 11 | 38 | | |
| | 2249 | 2 | 1 | 2 | 4 | 2 | 5 | 0 | 5 | 7 | 18 | 16 | 17 | 29 | 35 | 13 | 40 | X | |
| X | 2251 | 2 | 1 | 1 | 1 | 4 | 7 | 2 | 7 | 9 | 20 | 18 | 19 | 31 | 37 | 15 | 42 | X | X |
| | 2253 | 2 | 1 | 0 | 3 | 6 | 9 | 4 | 9 | 11 | 22 | 20 | 21 | 33 | 39 | 17 | 44 | | |
| | 2255 | 2 | 1 | 2 | 0 | 1 | 0 | 6 | 11 | 13 | 1 | 22 | 23 | 35 | 0 | 19 | 46 | | |

Table 6

Therefore, let P equal the set of core seed primes for a primorial M. Then by the combinatorial Multiplication Principle, we create a formula for computing the number of potential twin primes in primorial M (with $P_m$ equal to the largest core seed prime) by calculating the product of the primorial's core seed primes, excluding <u>two</u> residues (0 and $P_r$ -2) from each core seed prime, $P_r \in P$, in the calculation of the potential twin primes, and starting from 3 (the second prime factor, since the residue for odd numbers under Mod 2 is always equal to 1):

$$T(M) = [(3\text{-}2)(5\text{-}2)(7\text{-}2)\ldots(P_m \text{-}2)] = \prod_{Pr=3}^{Pm} (Pr \text{ -}2) \qquad (4)$$

The result of the formula is the number of potential conjecture solutions (in this case, the number of potential twin primes), T(M). It excludes the first three actual twin primes (5, 7 and 13) for primorial 11#

(since 3, 5, 7 and 11 have a zero residue in their modular signature under a core seed prime) and includes the number 1. Like the formula for potential primes in a primorial, the formula for the number of potential conjecture solutions (potential twin primes) is multiplicative because it is based on primorial stacking. That is, if M is a primorial whose prime factors are the set P = {2, 3, 5,…$P_M$} and N is the next larger primorial whose largest prime factor is $P_N$, then it follows that

$$T(N) = T(M) * (P_N - 2) \tag{5}$$

and we will have a relatively consistent distribution of the potential twin primes of M within the primorial N.

## Part 3: Application of Modular Signatures to the Goldbach Conjecture

Modular Residue Cycles
Consider the behavior of the residues of the integers under modular arithmetic. As we move down the column of integers, incrementing each by adding 1 to the preceding value, we move through the modular remainders in each modular signature in a synchronized way.

As Table 7 below illustrates, the odd integers have modular residues under the primes that follow this repeating sequence: odd residues for the prime, followed by 0, then the even residues for the prime. For example, for seed prime = 7, the odd integers 1-13 have this sequence of mod 7 residues: 1, 3, 5, 0, 2, 4, 6. In contrast, the even integers 2-14 have the same modular residues, but with the even and odd residues flipped, and 0 placed at the end of the sequence (i.e., the mod 7 residues for the even integers 2-14 have this sequence: 2, 4, 6, 1, 3, 5, 0). Furthermore, as the integers increase, the modular residues repeat in the same sequence for both the odd and even integers for each cycle equal to twice the seed prime (necessary to pass through all of the modular residues for both the odd and even integers).

|  | Odd Integers Mod Residues | | | | Even Integers Mod Residues | | | | |
|---|---|---|---|---|---|---|---|---|---|
| Integers | 2 | 3 | 5 | 7 | 2 | 3 | 5 | 7 | |
| 1 | 1 | 1 | 1 | 1 | | | | | <-- Signature |
| 2 | | | | | 0 | 2 | 2 | 2 | <-- Signature |
| 3 | 1 | 0 | 3 | 3 | | | | | |
| 4 | | | | | 0 | 1 | 4 | 4 | |
| 5 | 1 | 2 | 0 | 5 | | | | | |
| 6 | | | | | 0 | 0 | 1 | 6 | |
| 7 | 1 | 1 | 2 | 0 | | | | | |
| 8 | | | | | 0 | 2 | 3 | 1 | |
| 9 | 1 | 0 | 4 | 2 | | | | | |
| 10 | | | | | 0 | 1 | 0 | 3 | |
| 11 | 1 | 2 | 1 | 4 | | | | | |
| 12 | | | | | 0 | 0 | 2 | 5 | |
| 13 | 1 | 1 | 3 | 6 | | | | | |
| 14 | | | | | 0 | 2 | 4 | 0 | |
| 15 | 1 | 0 | 0 | 1 | | | | | |
| 16 | | | | | 0 | 1 | 1 | 2 | |
| 17 | 1 | 2 | 2 | 3 | | | | | |
| 18 | | | | | 0 | 0 | 3 | 4 | |
| 19 | 1 | 1 | 4 | 5 | | | | | |
| 20 | | | | | 0 | 2 | 0 | 6 | |
| 21 | 1 | 0 | 1 | 0 | | | | | |
| 22 | | | | | 0 | 1 | 2 | 1 | |
| 23 | 1 | 2 | 3 | 2 | | | | | |
| 24 | | | | | 0 | 0 | 4 | 3 | |
| 25 | 1 | 1 | 0 | 4 | | | | | |
| 26 | | | | | 0 | 2 | 1 | 5 | |
| 27 | 1 | 0 | 2 | 6 | | | | | |
| 28 | | | | | 0 | 1 | 3 | 0 | |

Table 7

For Goldbach, of critical importance is the relationship of the residues of smaller odd integers to the residues of a larger even integer. For example, consider the odd integers less than 14 shown in Table 7. Because of the way the modular signatures work, all of the residues under Mod 7 for the smaller odd integers are less than or equal to the integer 14's Mod 7 residue (0, which is the same as 7 in the logic of modular arithmetic).

Further, because of the repeated cycling of the modular residues for 3, 5 and 7, several of the odd integers less than 14 will be prime in addition to the seed primes (3, 5 and 7). Later, when we need to find odd primes less than an even integer like 14, we can see here why they will exist, whatever even integer we pick – because the repeated cycling of the modular residues for the seed primes (core and non-core), offset by the differing gaps between the seed primes, will produce modular signatures whose residues will all be non-zero for some odd integers in each cycle of the largest seed prime through the odd integers.

Table 7 illustrates why at least one odd integer O1 always exists less than an even integer E, with O1 and E in different congruence classes for all primes less than $\sqrt{E}$. For Goldbach, we need to know if at least one odd <u>prime</u> number P1 < E exists always, with P1 and E in different congruence classes for all seed primes less than $\sqrt{E}$. There are three reasons why that will be true.

First, each seed prime operates a repeating cycle consisting of its residues through successive integers, all starting where N = 1. As we increment N by 1 (i.e. cycle through the integers), each mod residue increases by 1. And each seed prime, $P_r$ in turn, will complete a full cycle (equal in length to twice the prime for the odd integers) at a different point, because each possesses $P_r$ residues, with the difference in the ending integer between two primes determined by the size of the gap between the two primes. Further, every odd integer has a residue = 1 for mod 2 and every even integer has residue = 0 for Mod 2.

Second, as the integers cycle through each even integer, the smaller odd integers have already cycled through the smaller modular residues (i.e. congruence classes) for each seed prime (e.g. 2, 3, 5, 7, 11 and 13 for primorial 210). For example, for Mod 7, there are seven different congruence classes (0, 1, 2, 3, 4, 5, 6). Whichever one includes E, there are six others available for an odd prime, P1, and all six will occur in sequence in Mod cycles earlier than E, as shown in Table 8 below for Mod 7's cyclic residues.

| Integer | 1 | 2 | 3 | 4 | 5 | 6 | 7 | 8 | 9 | 10 | 11 | 12 | 13 | 14 |
|---|---|---|---|---|---|---|---|---|---|---|---|---|---|---|
| | | | | | | Mod 7 Residues | | | | | | | | |
| - Odd | 1 | | 3 | | 5 | | 0 | | 2 | | 4 | | 6 | |
| - Even | | 2 | | 4 | | 6 | | 1 | | 3 | | 5 | | 0 |

Table 8

Third, every even integer, E > 6, falls between two primorials: $2*3…P_N = E_L \leq E \leq 2*3…P_N*P_{N+1} = E_H$. Then the following are true about the modular signatures of all even integers in the lower bound, $E_L$:

1. The modular signatures of the even integers ≤ $E_L$ are all unique, using the modular residues for seed primes through the largest prime < $\sqrt{E_H}$, as proved in Theorem 1.
2. All odd integers < E cycle through the modular residues for each seed prime, $P_r$, in an order that differs from the order for even integers. Starting with N = 1, odd integers cycle through the sequence of residues 1,3,…$P_r$-2,0,2,4,…$P_r$-1, while the even integers cycle through the sequence of residues 2,4,…$P_r$-1, 1,3,…$P_r$-2,0. Because the residue class cycles through an odd integer before its next greater even integer, the residue of the odd integer immediately less than E will be different under each mod $P_r$ as shown in Tables 7 and 8. In fact, the odd and even integers in the range X to X + 2$P_r$ will only have the same modular residue once for each $P_r$ in the set of seed primes, leaving $P_r-1$ differences in the modular residues of an odd integer and an even integer in that range.

3. All of the primes $< E_L$ remain available to be potential Goldbach solutions for E. Those that fall below E/2 (including seed primes) can serve as a potential P1 in the pair P1 + P2 = E.

   For example, let E = even integer in this primorial range: 19# < E < 23#. Then all 646,029 primes less than 19# = 9,699,690 remain primes when computed mod 14929, the largest seed prime for 23# = 223,092,870, and their unique modular signatures are merely extended to include a non-zero modular residue from the sequence of mod 14929 prime residues for odd integers: 1, 3,…14929, 0, 2, 4,…14928. In fact, there will be $\lfloor$ 223,092,870 / 29,858 $\rfloor$ = 7,471 complete cycles of the 14,929 residues through the odd integers under mod 14929, with the 0-residue skipping all of the existing primes less than 9,699,690.

Therefore, because the product of all core seed primes of 19# = 9,699,690 (the lower bound in our example) will generate all possible combinations of the modular residues for the primes 2 through 19, and extending the modular signatures of the primes < 9,699,690 to include the non-zero residues of the primes through 14,929, we can expect (and must prove) that there will be at least one prime number P1 less than $E_L$ (= 9,699,690 in the example) whose modular signature will have prime residues that all differ respectively from those of E. Of course, this logic can extend to any even integer E by simply increasing the number of seed primes used for its modular signature under a larger primorial.

Using Congruence Classes
For our purposes, let E refer to the even integer that is being subjected to the test for Goldbach solution pairs, and let P1 and P2 refer to the two primes that are being tested, with P1 ≤ E/2 and P2 ≥ E/2. Under mod 3, there are three general cases for Goldbach pair solutions that sum to E > 4, as shown in Table 9 below:

| Mod 3 Solution Combinations | | |
|---|---|---|
| Even>4 | P1 | P2 |
| [0] | [1] | [2] |
|  | [2] | [1] |
| Exception: 6 = 3+3 | | |
| [1] | [2] | [2] |
|  | 3 | [1] |
| [2] | [1] | [1] |
|  | 3 | [2] |

Table 9

1. If an even integer E > 4 is in congruence class [0] Mod 3, then one member of the Goldbach pair would have to come from congruence class [1] and one from congruence class [2]. Of course, if E = 6, then 3+3 is the Goldbach solution pair (the lone exception for congruence class [0]).
2. Any even integer E > 4 in congruence class [1] Mod 3 would have to be the sum of 3 + a prime number from congruence class [1], or two primes from congruence class [2]. For example, 40 = 3 + 37, and 40 = 11 + 29.
3. Any even integer E > 4 in congruence class [2] Mod 3 would have to be the sum of 3 + a prime number from [2], or two primes from congruence class [1]. For example, 44 = 3 + 41, and 44 = 7 + 37.

Based on the foregoing, we can expect to find more Goldbach solutions for E > 4 in congruence class [0] Mod 3 than in [1] Mod 3 or [2] Mod 3 because all of the primes in [1] and [2] are available for potential matches as Goldbach pairs when E is in [0].

## Number of Prime Pair Solutions of Goldbach

Table 10 below lists the number of prime pairs that satisfy the Goldbach conjecture for E ≤ 210, a primorial (210 = 2*3*5*7). This small sample of the known solutions illustrates that the number of solutions varies for each successive even integer.

Number of Goldbach Pairs for 4 < E <= 210

| Even | Mod 3 | # Pairs | Even | Mod 3 | # Pairs | Even | Mod 3 | # Pairs | Even | Mod 3 | # Pairs |
|---|---|---|---|---|---|---|---|---|---|---|---|
| 6 | 0 | 1 | 58 | 1 | 4 | 110 | 2 | 6 | 162 | 0 | 10 |
| 8 | 2 | 1 | 60 | 0 | 6 | 112 | 1 | 7 | 164 | 2 | 5 |
| 10 | 1 | 2 | 62 | 2 | 3 | 114 | 0 | 10 | 166 | 1 | 6 |
| 12 | 0 | 1 | 64 | 1 | 5 | 116 | 2 | 6 | 168 | 0 | 13 |
| 14 | 2 | 2 | 66 | 0 | 6 | 118 | 1 | 6 | 170 | 2 | 9 |
| 16 | 1 | 2 | 68 | 2 | 2 | 120 | 0 | 12 | 172 | 1 | 6 |
| 18 | 0 | 2 | 70 | 1 | 5 | 122 | 2 | 4 | 174 | 0 | 11 |
| 20 | 2 | 2 | 72 | 0 | 6 | 124 | 1 | 5 | 176 | 2 | 7 |
| 22 | 1 | 3 | 74 | 2 | 5 | 126 | 0 | 10 | 178 | 1 | 7 |
| 24 | 0 | 3 | 76 | 1 | 5 | 128 | 2 | 3 | 180 | 0 | 14 |
| 26 | 2 | 3 | 78 | 0 | 7 | 130 | 1 | 7 | 182 | 2 | 6 |
| 28 | 1 | 2 | 80 | 2 | 4 | 132 | 0 | 9 | 184 | 1 | 8 |
| 30 | 0 | 3 | 82 | 1 | 5 | 134 | 2 | 6 | 186 | 0 | 13 |
| 32 | 2 | 2 | 84 | 0 | 8 | 136 | 1 | 5 | 188 | 2 | 5 |
| 34 | 1 | 4 | 86 | 2 | 5 | 138 | 0 | 8 | 190 | 1 | 8 |
| 36 | 0 | 4 | 88 | 1 | 4 | 140 | 2 | 7 | 192 | 0 | 11 |
| 38 | 2 | 2 | 90 | 0 | 9 | 142 | 1 | 8 | 194 | 2 | 7 |
| 40 | 1 | 3 | 92 | 2 | 4 | 144 | 0 | 11 | 196 | 1 | 9 |
| 42 | 0 | 4 | 94 | 1 | 5 | 146 | 2 | 6 | 198 | 0 | 13 |
| 44 | 2 | 3 | 96 | 0 | 7 | 148 | 1 | 5 | 200 | 2 | 8 |
| 46 | 1 | 4 | 98 | 2 | 3 | 150 | 0 | 12 | 202 | 1 | 9 |
| 48 | 0 | 5 | 100 | 1 | 6 | 152 | 2 | 4 | 204 | 0 | 14 |
| 50 | 2 | 4 | 102 | 0 | 8 | 154 | 1 | 8 | 206 | 2 | 7 |
| 52 | 1 | 3 | 104 | 2 | 5 | 156 | 0 | 11 | 208 | 1 | 7 |
| 54 | 0 | 5 | 106 | 1 | 6 | 158 | 2 | 5 | 210 | 0 | 19 |
| 56 | 2 | 3 | 108 | 0 | 8 | 160 | 1 | 8 | | | |

Table 10

Observe further that as we move through Table 10, the number of Goldbach solution pairs is <u>not</u> evenly distributed across the modular remainders for Mod 3: [0], [1] and [2], as illustrated in Figure 2 below.

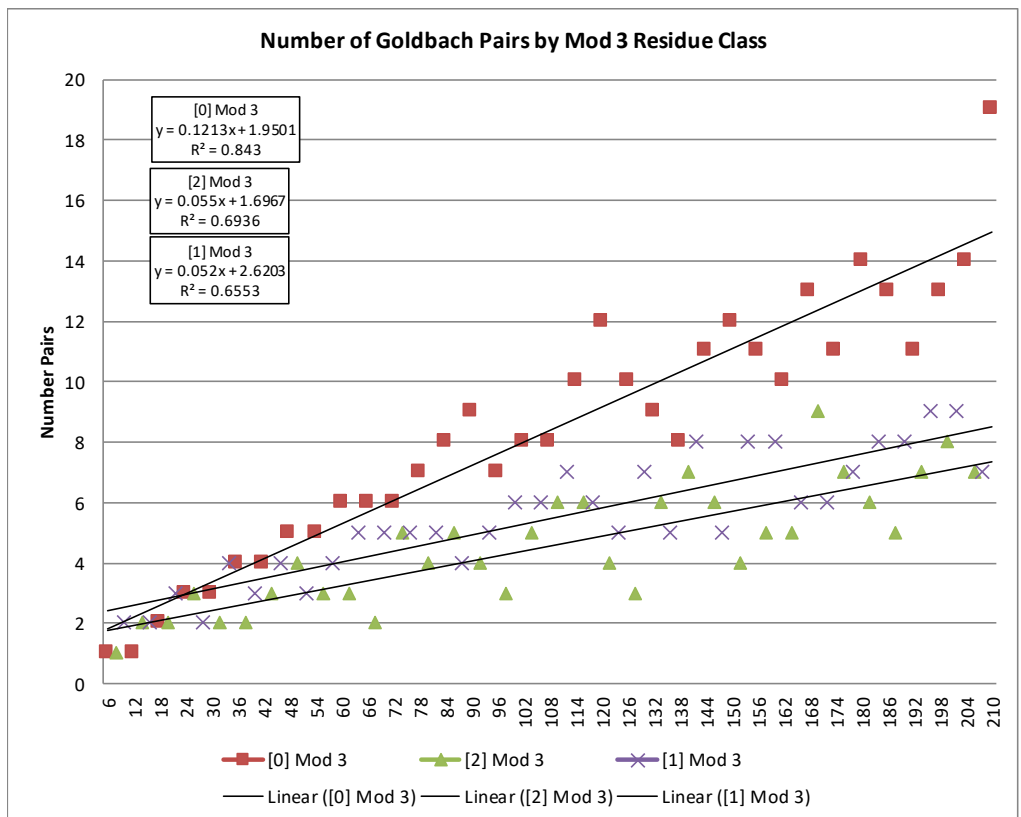

Figure 2

In fact, an even number in residue class [0] Mod 3 begins to generate more solution pairs than an even number in either [1] or [2] Mod 3. That result is entirely consistent with the solution rules contained in Table 9 because even numbers in [0] Mod 3 can draw on all odd primes in [1] and [2] Mod 3 for a solution, whereas even numbers in [1] or [2] Mod 3 can only draw on half of the available odd primes.

In Figure 2 we can see that the slope of the line of best fit for the number of Goldbach pairs in [0] Mod 3 is more than twice the slope of the line of best fit for the number of pairs in [1] Mod 3 or [2] Mod 3, whereas the slopes of the lines of best fit for [1] and [2] Mod 3 are fairly close, as expected.

Modular Residue Additive Combinations

If P is a prime number, then the set $Z_P = \{0, 1, …, P-1\}$ of integers contains all of the residues under mod P. For example, for P = 7, the modular residues equal the set {0, 1, 2, 3, 4, 5, 6}. Further, for any $x \in Z_P$, the additive inverse of x, -x, is also in $Z_P$ such that x + (-x) = 0. These properties of modular residues mean that for each prime number P, we can create a table showing all possible combinations of its modular residues that will sum to elements in $Z_P$.

| 7 | 0 | 1 | 2 | 3 | 4 | 5 | 6 |
|---|---|---|---|---|---|---|---|
| 0 | **0** | 1 | 2 | 3 | 4 | 5 | 6 |
| 6 | 6 | **0** | 1 | 2 | 3 | 4 | 5 |
| 5 | 5 | 6 | **0** | 1 | 2 | 3 | 4 |
| 4 | 4 | 5 | 6 | **0** | 1 | 2 | 3 |
| 3 | 3 | 4 | 5 | 6 | **0** | 1 | 2 |
| 2 | 2 | 3 | 4 | 5 | 6 | **0** | 1 |
| 1 | 1 | 2 | 3 | 4 | 5 | 6 | **0** |

Table 11

As illustrated in Table 11 for Mod 7, there are 49 combinations of the 7 residues of Mod 7, and 7 of those produce 0. In general, for any prime P, there will be $P^2 – P = P*(P - 1)$ combinations that will sum to a modular residue <>0.

The implication of Table 11 is that since an even integer E > 4 has a modular signature under primes 2 through P, it will be possible to find at least two combinations of modular residues for each P > 2 that will sum to the modular residue of E for each seed prime in E's modular signature under the smallest primorial ≥ E. Table 12 below illustrates this fact for E = 60, and P1 = 17 and P2 = 43. In the modular signature of 17 under seed primes 2, 3, 5, 7, 11 and 13, the modular residues of 17 (1, 2, 2, 3, 6, 4) all differ from the corresponding modular residues in the modular signature of E = 60 (0, 0, 0, 4, 5, 8).

|   |   |   | Mod Residues | | | | | | |
|---|---|---|---|---|---|---|---|---|---|
|   |   |   | 2 | 3 | 5 | 7 | 11 | 13 | <- Seed Primes |
| Even | = | 60 | 0 | 0 | 0 | 4 | 5 | 8 | <- Signature |
| - P1 | = | 17 | 1 | 2 | 2 | 3 | 6 | 4 |   |
| P2 | = | 43 | 1 | 1 | 3 | 1 | 10 | 4 |   |

Table 12

Therefore, a necessary condition for proving the Goldbach conjecture will be to show that for each even number E > 4, there must exist at least one pair of odd prime numbers less than E whose modular signatures under the smallest primorial ≥ E will add together to produce the modular signature of E under the primorial's seed primes, as illustrated in Table 12 for E = 60 = 17 + 43.

Restriction on Using Modular Arithmetic for Goldbach Matches

Table 9's decision rules for Mod 3 indicated that for each seed prime there are certain combinations of the Mod 3 modular residues that must be used, or must not be used, for either E – P1 = P2, or P1 + P2 = E, to ensure that both P1 and P2 are odd prime numbers whose sum will equal E.

Assume E > 4 is <u>not</u> equal to twice a prime (the trivial case for the Goldbach conjecture).
- If a modular residue of E equals 0 under seed prime $P_S$, with P1 < E/2, and <u>if P1 equals $P_S$</u>, then P2 (which is > E/2) will not be a prime number because the modular residue of P2 under seed prime $P_S$ would also have to equal 0.
- If a modular residue of E does <u>not</u> equal 0 under seed prime $P_S$, then P1 cannot have the same modular residue (i.e. be in the same congruence class) under seed prime $P_S$. Otherwise, when we compute Mod(E) – Mod(P1) = Mod (P2) under $P_S$, it will force the modular residue of P2 under $P_S$ to equal 0, and then P2 > E/2 will not be a prime number because it will be a composite of $P_S$.

These two conditions can be combined into a single necessary condition for finding a Goldbach solution: **Where even number E > 4 is not equal to twice a prime and P1 is a prime less than E/2, every respective residue in the modular signatures of P1 and E must be in a different congruence class of the associated seed primes for the smallest primorial ≥ E.** Note that the same requirement does not apply to P2, as it can have a residue that equals a corresponding residue of E if P1 is a seed prime.

Note too in Table 10's list of the number of Goldbach solutions for even integers up to 210 that the even integer 68 in [2] Mod 3 generated the smallest number of prime pair solutions (2) for even numbers > 38. Therefore, it will be instructive to examine what circumstances led to so few solutions.

Analysis of the Goldbach pairs indicates that there are only two prime pairs that sum to 68: 7 + 61 and 31 + 37, where the seed primes for 68's modular signature are 2, 3, 5, 7, 11 and 13 (because 30 < 68 < 210 and 13 is the largest prime < $\sqrt{210}$ ). Examine these two modular residue tables (similar to Table 12) to see why:

|  |  |  | Mod Residues | | | | | |  |
|---|---|---|---|---|---|---|---|---|---|
|  |  |  | 2 | 3 | 5 | 7 | 11 | 13 | <- Seed Primes |
| Even | = | 68 | 0 | 2 | 3 | 5 | 2 | 3 | <- Signature |
| - P1 | = | 7 | 1 | 1 | 2 | 0 | 7 | 7 |  |
| P2 | = | 61 | 1 | 1 | 1 | 5 | 6 | 9 |  |

Table 13

In Table 13, we see that the modular residues of E and P1 differ for each seed prime. Consequently, 7 and 61 satisfy the necessary conditions for the sum of a pair of primes to equal an even integer. Note, too, that the modular residue of P1 = 0 for mod 7. Since E = 68 is not a multiple of 7, that means that 68 and 7 have different residues under mod 7. So, a Goldbach pair can contain a seed prime so long as the even number is not a multiple of that seed prime, with one trivial exception (E = 2P, in which case P will be a Goldbach solution). Next, consider an alternate pair of primes for E in Table 14.

|  |  |  | Mod Residues | | | | | |  |
|---|---|---|---|---|---|---|---|---|---|
|  |  |  | 2 | 3 | 5 | 7 | 11 | 13 | <- Seed Primes |
| Even | = | 68 | 0 | 2 | 3 | 5 | 2 | 3 | <- Signature |
| - P1 | = | 31 | 1 | 1 | 1 | 3 | 9 | 5 |  |
| P2 | = | 37 | 1 | 1 | 2 | 2 | 4 | 11 |  |

Table 14

In Table 14, we see an alternate pair of primes which satisfy the Goldbach conjecture for E = 68, and in which neither one of the primes is a seed prime (i.e. both are greater than 13) under primorial 7# = 210. In this case, P1 = 31 is not in the same congruence class (i.e. has a different modular residue) as E = 68 for any of the seed primes: 2, 3, 5, 7, 11 and 13. Because all of the modular residues of 31 are not equal to 0 and it is not a seed prime, we know by Theorem 2 that it is a prime number. Further, we can reach the same conclusion for 37.

Now consider a case in which P1 is prime, but P2 is not a prime because the choice of P1 violates the necessary condition for a Goldbach pair.

|  |  |  | Mod Residues |  |  |  |  |  |  |
|---|---|---|---|---|---|---|---|---|---|
|  |  |  | 2 | 3 | 5 | 7 | 11 | 13 | <- Seed Primes |
| Even | = | 68 | 0 | 2 | 3 | 5 | 2 | 3 | <- Signature |
| - P1 | = | 19 | 1 | 1 | 4 | 5 | 8 | 6 |  |
| P2 | = | 49 | 1 | 1 | 4 | 0 | 5 | 10 |  |

Table 15

In Table 15, we can see that the Mod 7 residue for both E = 68 and P1 = 19 equals 5 (i.e. they are in the same congruence class [5] for Mod 7). That forces P2 to have a Mod 7 residue = 0, which means that P2 is a multiple of 7 and therefore is not prime because it is divisible by 1, 7 and 49.

Since 68 is in congruence class [2] under Mod 3, we know that its Goldbach solutions have to come from either 3 plus a prime in congruence class [2] Mod 3 or from two primes in congruence class [1] Mod 3. However, since 3 has the same modular residue as 68 under Mod 5, we know that 3 cannot be a member of the Goldbach solution set. Therefore, any Goldbach solution for E = 68 must come from congruence class [1] Mod 3 – and both solution pairs do: (7+61) and (31+37).

Therefore, for each seed prime, there are <u>at most</u> two congruence classes that will eliminate the odd number from being a prime number that can be a Goldbach solution in a cycle of the seed prime's residues:
- A zero residue means the odd integer is a composite if the odd integer is not a seed prime.
- A non-zero residue that matches the corresponding seed prime residue in E which we are testing for Goldbach solutions.

Using Combinations To Count Residue Matches For Goldbach Solutions
Since we will be eliminating up to two potential residue matches from each congruence class of the modular signature of an even number E, we can use the same combinatorial formula, Equation (4), to count the number of potential Goldbach solutions in a primorial M (with $P_M$ equal to the largest core seed prime) as we used to count potential twin prime solutions.

While the result of the formula in Equation (4) includes the number 1 (which has no 0 residues), it understates the actual number of potential Goldbach solutions for two reasons:
- The formula excludes the core seed primes since they are not potential primes (they each have a zero residue in their modular signature) – and core seed primes can be part of a Goldbach solution.
- The formula always excludes two residues for each core seed prime even though the core seed prime may be a factor of the even number (in which case we would only have had to eliminate one residue, 0, to get the potential Goldbach solution for that seed prime).

As was the case for potential twin primes, the formula for the number of potential Goldbach solutions is multiplicative because it is based on primorial stacking. That is, if M is a primorial whose prime factors are the set P = {2, 3, 5,…$P_M$} and N is the next larger primorial whose largest prime factor is $P_N$, then it follows that Equation (4) will provide us the required count of potential Goldbach solutions in a primorial. Similarly, because T(N) = T(M) * ($P_N$ – 2), we will have a relatively even distribution of the potential Goldbach solutions of M within the primorial N, as the set of potential Goldbach solutions is a subset of the set of potential primes, as was the case for potential twin primes.

The foregoing suggests that it should be possible to identify the potential primes that will have no residue matches with the applicable core seed primes and therefore be potential Goldbach solutions, since, for any seed prime, at most two residues are eliminated: 0 and the matching residue in the modular signature of E for each seed prime of the smallest primorial greater than or equal to E.

**Theorem 5**: Let A equal a primorial whose largest prime factor (and therefore, its largest core seed prime) is $P_N$. If $P_Z$ equals the smallest <u>non-core</u> seed prime of A, the potential primes (and their subsets) that are less than $(P_Z)^2$ in A will not be converted to composites; they will remain prime numbers.

Proof.
1. Let Q = the set of seed primes for A, with $P_r \in Q$.
2. Under the combinatorial Multiplication Principle, the core seed primes of A, (2, 3, 5,…$P_N$) will generate all possible combinations of modular residues for the integers from 1 to A under its core seed primes. The number of potential prime numbers in A will equal $\varphi(A)$ and by <u>Equation (4)</u> the number of potential twin primes or Goldbach solutions in A will equal T(A) = $\prod_{Pr=3}^{Pn}$ (Pr-2) for primes Pr = 3 to $P_N$.
3. There will be no zero residues in the modular signatures of all potential primes (or the elements in their subsets) by virtue of the formula used to generate them under the core seed primes in A.
4. By <u>Theorem 3</u>, since $P_Z$ equals the smallest non-core seed prime of A, $(P_Z)^2$ will be the first new composite generated by the non-core seed primes of A.
5. Since $(P_Z)^2$ will be the first new composite generated by the non-core seed primes of A, by <u>Corollary 3.2</u> the modular signatures of all potential primes (including potential primes which are potential twin primes or potential Goldbach solutions) that are less than $(P_Z)^2$ will not have any zero modular residues under the non-core seed primes of A; instead, the modular signatures of all potential primes and potential twin primes (or potential Goldbach solutions) that are less than $(P_Z)^2$ will only contain non-zero residues under the non-core seed primes of A.
6. Since any potential primes and potential twin primes (or Goldbach solutions) < $(P_Z)^2$ in A will have no zero residues in their modular signatures, by <u>Theorem 2</u> they will remain prime numbers.

## Part 4: Application of Primorial Stacking to Conjecture Solutions

Next, consider how the number of potential twin primes or potential Goldbach solutions generated by Equation (4) in each primorial changes under primorial stacking. That is, when we consider two or more primorials, will the average number of potential twin primes or Goldbach solutions in the smaller primorial remain constant in each cycle within the larger primorial, or be smaller because of the additional core seed primes in the larger primorial?

Table 16 below illustrates the increasing number of potential twin prime / Goldbach solutions in each primorial as the primorials grow. The number of potential primes is calculated using Euler's Totient Formula. The number of potential twin primes / Goldbach solutions is calculated using Equation (4), with each potential twin or Goldbach solution factor in column 6 being 2 less than the corresponding prime factor in column 2 for the primorial in column 3.

| Count | Prime P | Primorial P# | Number Of Cycles By Prior Primorial | Potential Prime Factors P-1 | Potential Equation (4) Solution Factors P-2 | Number Potential Primes In Primorial | Number Potential Solutions In Primorial | Avg Number Potential Solutions In Each Prior Primorial Cycle |
|---|---|---|---|---|---|---|---|---|
| 1 | 2 | | | | | | | |
| 2 | 3 | 6 | | 2 | 1 | | | |
| 3 | 5 | 30 | | 4 | 3 | 8 | 3 | |
| 4 | 7 | 210 | 7 | 6 | 5 | 48 | 15 | 2.1 |
| 5 | 11 | 2,310 | 11 | 10 | 9 | 480 | 135 | 12.3 |
| 6 | 13 | 30,030 | 13 | 12 | 11 | 5,760 | 1,485 | 114.2 |
| 7 | 17 | 510,510 | 17 | 16 | 15 | 92,160 | 22,275 | 1,310.3 |
| 8 | 19 | 9,699,690 | 19 | 18 | 17 | 1,658,880 | 378,675 | 19,930.3 |
| 9 | 23 | 223,092,870 | 23 | 22 | 21 | 36,495,360 | 7,952,175 | 345,746.7 |
| 10 | 29 | 6,469,693,230 | 29 | 28 | 27 | 1,021,870,080 | 214,708,725 | 7,403,749.1 |

Table 16

Observe the following in the table:
- Column 4 shows the number of cycles of the prior primorial in the current primorial. That is, there are 7 cycles of primorial 30 in primorial 210 since 210 = 30*7.
- Column 8 shows that the number of potential twin prime or Goldbach solutions increases as each primorial increases. Column 9 divides the number of potential solutions by the number of cycles of the prior primorial (column 4). That calculation shows that the average number of potential twin primes or Goldbach solutions in each cycle of the prior primorial is less than the number of the potential solutions in the prior primorial. In other words, if we look at the effect of primorial stacking on the number of potential twin primes or Goldbach solutions, the addition of the extra core seed prime in each new primorial reduces the number of potential twin primes or Goldbach solutions in each prior primorial as it cycles within the larger primorial's more expansive modular signatures.

For example, the number of potential solutions in primorial 2310 is = 135. Primorial 2310 then has 13 cycles in primorial 30030. If each of those cycles contained 135 potential solutions, then there would be 13 * 135 = 1755 potential solutions in primorial 30030. However, there are only 1485 potential solutions in primorial 30030. The reduction from 135 is expected in the cycles of primorial 2310 within primorial 30030, since some of the modular signatures of the potential twin primes or Goldbach solutions from primorial 2310 will receive a zero residue from the additional core seed prime in 30030. Therefore, the average number of potential solutions in M = 2310 is reduced to 135*([13-2]/13) = 135*(11/13) = 114 in each cycle of the 13 cycles of M within N = 30030 so that the total number of potential twin primes or Goldbach solutions in N will remain T(M) * (13-2) = 1485.

While the average number of potential twin primes or Goldbach solutions in a primorial is reduced within a larger primorial, note that there may or may not be potential solutions in every cycle, with only a small variance between each, as shown in Table 5, depending on the relative size of the two primorials. Theorem 6 and its lemmas addresses how to ensure that the average number of potential twin primes or Goldbach solutions in a primorial is greater than 1 as it cycles within a larger primorial.

**Theorem 6**: Let M equal a primorial > 30 and let N equal a primorial greater than M. Let $P_M$ equal the largest prime factor of M, let $P_S$ be the smallest prime number greater than $P_M$ and let $P_Z$ equal the largest prime factor of N (i.e., N = M*$P_S$…$P_Z$). Let T($M_M$) equal the number of potential conjecture solutions in M under the core seed primes of M, and let T($M_N$) equal the average number of potential solutions in each cycle of M under the core seed primes of N. Then, the number of potential solutions of N equals T(N) = T($M_M$) * $\prod_{P_S}^{P_Z}(Pr - 2)$, where Pr equals a prime number from $P_S$ to $P_Z$, and <u>on average</u> there will be T($M_N$) = T(N) / $\prod_{P_S}^{P_Z}(Pr)$ = T($M_M$) * $\prod_{P_S}^{P_Z}((Pr - 2)/Pr)$ potential conjecture solutions in each cycle of M in N.

Proof.
1. N = 2*3*5*…$P_M$*$P_S$…$P_Z$.
2. By Equation (4), T(N) = $\prod_3^{Pz}$ (Pr − 2) for Pr = primes 3 to $P_Z$
3. By Equation (4), T($M_M$) = $\prod_3^{Pm}$(Pr − 2) for Pr = primes 3 to $P_M$
4. T(N) = $\prod_3^{Pm}(Pr-2)$ * $\prod_{Ps}^{Pz}(Pr-2)$
5. T(N) = T($M_M$) * $\prod_{Ps}^{Pz}(Pr-2)$ over the core seed primes of N.
6. Since there are $\prod_{Ps}^{Pz}(Pr)$ cycles of M in N, on average there will be T($M_N$) = T(N) / $\prod_{Ps}^{Pz}(Pr)$ potential solutions in each cycle of M in N.
7. T($M_N$) = T(N) / $\prod_{Ps}^{Pz}(Pr)$
8. T($M_N$) = T($M_M$) * $\prod_{Ps}^{Pz}(Pr-2)$ / $\prod_{Ps}^{Pz}(Pr)$
9. T($M_N$) = T($M_M$) * $\prod_{Ps}^{Pz}((Pr-2)/Pr)$ potential solutions in each cycle of M in N on average.

Managing the Product Factor
In Theorem 6, the product factor of T($M_N$) is equal to $\prod_{Ps}^{Pz}((Pr-2)/Pr)$ and is less than 1. It is used to compute the average number of potential twin primes or Goldbach solutions of a smaller primorial, M, in a larger primorial, N. A careful reading of the formula would show that if the reciprocal of the product factor is greater than T($M_M$), then the average number of potential twin primes of M in N, T($M_N$), will be less than 1. Therefore, it will be necessary to ensure that $P_Z$ (the largest core seed prime of primorial N), is small enough to enable the reciprocal of the product factor to be less than T($M_M$).

**Lemma 6.1**: Let M equal a primorial > 30 and let N equal a primorial greater than M. Let $P_M$ equal the largest prime factor of M, let $P_S$ be the smallest prime number greater than $P_M$ and let $P_Z$ equal the largest prime factor of N (i.e., N = M*$P_S$…$P_Z$), <u>where $P_Z$ is the largest prime number less than $\sqrt{M}$</u>. Let T($M_M$) equal the number of potential solutions in M using Equation (4) under the core seed primes of M, and let T($M_N$) equal the average number of potential conjecture solutions in each cycle of M under the core seed primes of N. Then T($M_M$) will be greater than the reciprocal of the product factor, $\prod_{Ps}^{Pz}((Pr-2)/Pr)$.

Proof.
1. T($M_M$) > $\sqrt{M}$ / $P_M$
2. By definition, $\sqrt{M}$ > $P_Z$.
3. $\sqrt{M}$ / $P_M$ > $P_Z$ / $P_M$.
4. $P_Z$ / $P_M$ > $\prod_{Ps}^{Pz}(Pr / (Pr-2)$ [14]
5. T($M_M$) > $\prod_{Ps}^{Pz}(Pr / (Pr-2))$

Table 17 below illustrates the relationship of the various factors used in Lemma 6.1, showing how using the square root of M to limit the size of $P_Z$ (and N) ensures that the reciprocal of the product factor will be less than T($M_M$).

| | Primorials | | Largest Core Seed Prime | | Potential Equation (4) Solutions of M | | | Reciprocal of Product Factor |
|---|---|---|---|---|---|---|---|---|
| | M | N | Pm | Pz | T(M) in M | Product Factor | Avg T(M) in N | |
| 1 | 210 | 13# | 7 | 13 | 15 * | 0.692308 = | 10 | 1.44 |
| 2 | 2,310 | 47# | 11 | 47 | 135 * | 0.436373 = | 59 | 2.29 |
| 3 | 30,030 | 173# | 13 | 173 | 1,485 * | 0.307356 = | 456 | 3.25 |
| 4 | 510,510 | 709# | 17 | 709 | 22,275 * | 0.218553 = | 4,868 | 4.58 |
| 5 | 9,699,690 | 3,109# | 19 | 3,109 | 378,675 * | 0.164156 = | 62,162 | 6.09 |
| 6 | 223,092,870 | 14,929# | 23 | 14,929 | 7,952,175 * | 0.126197 = | 1,003,543 | 7.92 |
| 7 | 6,469,693,230 | 80,429# | 29 | 80,429 | 214,708,725 * | 0.098251 = | 21,095,426 | 10.18 |
| 8 | 200,560,490,130 | 447,829# | 31 | 447,829 | 6,226,553,025 * | 0.079161 = | 492,902,698 | 12.63 |
| 9 | 7,420,738,134,810 | 2,724,079# | 37 | 2,724,079 | 217,929,355,875 * | 0.064543 = | 14,065,843,393 | 15.49 |

Table 17

Therefore, if the largest core seed prime (prime factor) of primorial N is less than the square root of primorial M, then we can be assured that the average number of potential twin primes or Goldbach solutions of M in N is greater than 1.

Further, we can show that a given T(M$_N$) will be larger than its predecessors, as is apparent in Table 17 above and Table 18 below, as shown by the ratio of successive values of T(M) in N in column 6.

| | Primorials | | Avg Equation (4) Potential Solutions | | Ratio of Successive T(M) | | |
|---|---|---|---|---|---|---|---|
| | M | N | T(M1) in N1 | T(M2) in N2 | T(M2) in N2 / T(M1) in N1 | T(M2) in M2 / T(M1) in M1 | PF2 / PF1 |
| 1 | 210 | 13# | | | | | |
| 2 | 2,310 | 47# | 10 | 59 | 5.6728 = | 9 * | 0.6303 |
| 3 | 30,030 | 173# | 59 | 456 | 7.7478 = | 11 * | 0.7043 |
| 4 | 510,510 | 709# | 456 | 4,868 | 10.6661 = | 15 * | 0.7111 |
| 5 | 9,699,690 | 3,109# | 4,868 | 62,162 | 12.7688 = | 17 * | 0.7511 |
| 6 | 223,092,870 | 14,929# | 62,162 | 1,003,543 | 16.1441 = | 21 * | 0.7688 |
| 7 | 6,469,693,230 | 80,429# | 1,003,543 | 21,095,426 | 21.0209 = | 27 * | 0.7786 |
| 8 | 200,560,490,130 | 447,829# | 21,095,426 | 492,902,698 | 23.3654 = | 29 * | 0.8057 |
| 9 | 7,420,738,134,810 | 2,724,079# | 492,902,698 | 14,065,843,393 | 28.5368 = | 35 * | 0.8153 |

Table 18

In the table, the first ratio uses M1 = 210 and M2 = 2310, the second ratio uses M1 = 2310 and M2 = 30,030, etc. Note that as the primorials grow, the ratio of their respective product factors (shown in the last column) also increases. That result guarantees that the ratio of the successive T(M) in N will increase.

**Lemma 6.2**:
Let M1 equal a primorial > 30 and let N1 equal a primorial greater than M1.
- Let P$_{Z1}$ equal the largest prime factor of N1, where P$_{Z1}$ is the largest prime number less than $\sqrt{M1}$.
- Let P$_{M1}$ equal the largest prime factor of M1.
- Let P$_{S1}$ be the smallest prime number greater than P$_{M1}$ (i.e., N1 = M1*P$_{S1}$…P$_{Z1}$).

Let M2 equal the next primorial greater than M1 and let N2 equal a primorial greater than M2.
- Let P$_{Z2}$ equal the largest prime factor of N2, where P$_{Z2}$ is the largest prime number less than $\sqrt{M2}$.
- Let P$_{M2}$ equal the largest prime factor of M2. P$_{M2}$ = P$_{S1}$, and M2 = M1 * P$_{S1}$.
- Let P$_{S2}$ be the smallest prime number greater than P$_{M2}$ (i.e., N2 = M2*P$_{S2}$…P$_{Z2}$).

Let T(M1$_{M1}$) equal the number of potential conjecture solutions in M1 using Equation (4) under the core seed primes of M1, and let T(M1$_{N1}$) equal the average number of potential solutions in each cycle of M1 under the core seed primes of N1. Similarly, let T(M2$_{M2}$) equal the number of potential solutions in M2 under the core seed primes of M2, and let T(M2$_{N2}$) equal the average number of potential solutions in each cycle of M2 under the core seed primes of N2. Then T(M2$_{N2}$) > T(M1$_{N1}$) and each successive T(M$_N$) will be greater than its predecessors.

Proof.
1. Since M2 = M1 * P$_{S1}$, T(M2) = T(M1) * (P$_{S1}$ – 2).
2. Let P$_{Z1+1}$ equal the smallest prime greater than P$_{Z1}$.
3. Then the ratio of T(M2$_{N2}$) to T(M1$_{N1}$) equals

$$\frac{[T(M1) * (P_{S1} - 2)] * \prod_{Ps2}^{Pz2}((Pr - 2)/Pr)}{T(M1) * \prod_{Ps1}^{Pz1}((Pr - 2)/Pr)} = P_{S1} * \prod_{Pz1+1}^{Pz2}((Pr - 2)/Pr)$$

5. P$_{S1}$ > $\sqrt{M2}$ / P$_{Z1}$ [15]

6. $\sqrt{M2} > P_{Z2}$
7. $\sqrt{M2} / P_{Z1} > P_{Z2} / P_{Z1}$
8. $P_{Z2} / P_{Z1} > \prod_{Pz1+1}^{Pz2} (Pr / (Pr - 2))$
9. $Ps_1 > \prod_{Pz1+1}^{Pz2} (Pr / (Pr - 2))$
10. Since $Ps_1$ is greater than the reciprocal of the product factor, the ratio of $T(M2_{N2})$ to $T(M1_{N1})$ is greater than one, and $T(M2_{N2}) > T(M1_{N1})$ and each successive $T(M_N)$ will be larger than its predecessors.

**Lemma 6.3**: Let M equal a primorial > 30 and let N equal a primorial greater than M. Let $P_M$ equal the largest prime factor of M, let $P_S$ be the smallest prime number greater than $P_M$ and let $P_Z$ equal the largest prime factor of N (i.e., $N = M*P_S…P_Z$), where $P_Z$ is the largest prime number less than $\sqrt{M}$. Let $T(M_M)$ equal the number of potential solutions (twin prime or Goldbach) in M using Equation (4) under the core seed primes of M, and let $T(M_N)$ equal the average number of potential solutions in each cycle of M under the core seed primes of N. Then $T(M_N) > 10$.

Proof.
1. By Equation (4), $T(M_M) = (3-2)(5-2)…(P_M -2)$
2. By Equation (4), $T(M_N) = (3-2)(5-2)…(P_M -2)(P_S - 2) …(P_Z - 2) = T(M_M) * \prod_{PS}^{Pz}(Pr - 2)$ where $P_S \le Pr \le P_Z$ and Pr is prime
3. $T(M_N) = 10$ for the combination of primorials (M = 210 and N = 30030), where 13 is the smallest prime < 210 and 13# = 30030.
4. Since by Lemma 6.2 every $T(M_N)$ is larger than its predecessor, $T(M_N) > 10$ for all $M \ge 210$, where the largest core seed prime of N is equal to the largest prime $< \sqrt{M}$.

Building the Primorial Scaffold for Twin Primes

Table 19 below shows groups of three primorials with a specific relationship. A is the smallest primorial, B is the next largest, and C is the primorial whose largest factor is the largest prime $< \sqrt{B}$. The largest core seed prime of A is labeled Pa; the largest of B is Pb, and the largest of C is Pc.

| | Primorials | | | Largest Core Seed Prime | | | Potential Twin Primes of A | | |
|---|---|---|---|---|---|---|---|---|---|
| | A | B | C | Pa | Pb | Pc | T(A) in A | Product Factor | Avg T(A) in C |
| 1 | 210 | 2,310 | 47# | 7 | 11 | 47 | 15 | 0.357032 | 5 |
| 2 | 2,310 | 30,030 | 173# | 11 | 13 | 173 | 135 | 0.260070 | 35 |
| 3 | 30,030 | 510,510 | 709# | 13 | 17 | 709 | 1,485 | 0.192841 | 286 |
| 4 | 510,510 | 9,699,690 | 3,109# | 17 | 19 | 3,109 | 22,275 | 0.146876 | 3,272 |
| 5 | 9,699,690 | 223,092,870 | 14,929# | 19 | 23 | 14,929 | 378,675 | 0.115224 | 43,632 |
| 6 | 223,092,870 | 6,469,693,230 | 80,429# | 23 | 29 | 80,429 | 7,952,175 | 0.091475 | 727,428 |
| 7 | 6,469,693,230 | 200,560,490,130 | 447,829# | 29 | 31 | 447,829 | 214,708,725 | 0.074054 | 15,900,087 |
| 8 | 200,560,490,130 | 7,420,738,134,810 | 2,724,079# | 31 | 37 | 2,724,079 | 6,226,553,025 | 0.061054 | 380,157,930 |

Table 19

Observe the following in the table:
- Starting with the second iteration of primorials A-B-C, each new primorial A equals the prior primorial B, with the largest prime $< \sqrt{B}$ equaling the largest core seed prime of primorial C.
- The number of potential twin primes in A, T(A), is calculated using Equation (4) and is shown in column 8. Those values are multiplied by the product factor in column 9 to get the average number of potential twin primes in each cycle of A within C, Avg T(A) in C in column 10. The product factor is equal to $\prod_{Pb}^{Pc}((Pr - 2)/Pr)$, where Pb is the largest core seed prime (factor) of B and Pc is the largest core seed prime (factor) of C. That is, the product factor is calculated using the prime factors that are greater than those in A and less than or equal to those in C.

Alternatively, the average value of T(A) in C can be derived by using Theorem 6: calculate T(C) and divide T(C) by the product of the core seed primes between Pb and Pc inclusive.
- Note that while the primorials A, B and C grow very rapidly, the largest core seed primes increase much more slowly. The consequence is that the product factor, which is less than 1 and is derived from the core seed primes, does not shrink as quickly as the primorials grow, so that the average T(A) in C continues to grow rapidly.

The product factor used to compute the average number of potential twin primes of primorial A in primorial C is equal to $\prod_{Pb}^{Pc}((Pr-2)/Pr)$, where Pb is the largest core seed prime of primorial B and Pc is the largest core seed prime of primorial C. For example, consider A = 17# = 510,510 and B = 19# = 9,699,690 and C = 3109#. The product factor for T(A) in C starts with Pb = 19 and grows until Pc = 3109. The average T(A) in C = T(A$_A$) * Product Factor = 22,275 * 0.14687644 = 3,272.

Table 19 illustrates again that, consistent with Lemma 6.2, each successive T(A) in C is larger than its predecessor under the specific construct of primorials A, B and C, where B is the smallest primorial greater than A, and the largest core seed prime (prime factor) of C is the largest prime $< \sqrt{B}$ (and is also the largest non-core seed prime of B).

**Theorem 7**: Let A > 30 be a primorial and let B be the smallest primorial larger than A. Let Pb = the largest core seed prime of B, such that B = A * Pb. Let P$_C$ be the largest prime $< \sqrt{B}$ and let primorial C = P$_C$#. Let T(A$_A$) equal the number of potential twin primes of A in A and T(A$_C$) equal the average number of potential twin primes in each cycle of A in C, and T(B$_C$) equal the average number of potential twin primes in each cycle of B in C. Then T(B$_C$) > T(A$_C$) and T(B$_C$) – T(A$_C$) > 5.

Proof.
1. Let Pa = the largest core seed prime of A.
2. By Equation (4), T(A$_A$) = (3-2)(5-2)…(Pa -2)
3. By Equation (4), T(B$_B$) = (3-2)(5-2)…(Pa -2)(Pb – 2) = T(A$_A$) * (Pb – 2)
4. By Theorem 6, T(A$_C$) = T(A$_A$) * $\prod_{Pb}^{Pc}((Pr-2)/Pr)$, where Pr = a prime number from Pb to Pc
5. By Theorem 6, T(B$_C$) = T(C) / $\prod_{Ps}^{Pc}(Pr)$, where Ps is the next larger prime than Pb
6. By Theorem 6, T(A$_C$) = T(C) / $\prod_{Pb}^{Pc}(Pr)$
7. Since Pb*$\prod_{Ps}^{Pc}(Pr)$ = $\prod_{Pb}^{Pc}(Pr)$, T(B$_C$) = T(A$_C$) * Pb
8. Since T(B$_C$) = Pb * T(A$_C$), then T(B$_C$) – T(A$_C$) = T(A$_C$)*Pb – T(A$_C$) = T(A$_C$) * (Pb – 1)
9. T(A$_C$) = 5 for the combination of primorials (A = 210, B = 2310 and C = 30030), where 210 is the smallest primorial > 30
10. Since T(B$_C$) = Pb * T(A$_C$), then T(B$_C$) – T(A$_C$) > 5 for A ≥ 210 since by Lemma 6.2 each successive Pb, T(B$_C$) and T(A$_C$) will be larger than their predecessors.

Theorems 6 and 7 provide key elements for the proof of the existence of an infinite number of twin primes because they demonstrate that there are potential twin primes in successive primorials, and where those potential twin primes are less than the square root of the largest non-core seed prime of a primorial, they are not converted to composites. By Theorem 5, they will remain twin primes. This is illustrated in Table 20 below, an extension of Table 19.

Table 20 illustrates that the average T(A) in C is less than the average T(B) in C, as proved in Theorem 7. Further, since T(B) > T(A) under Equation (4), there are potential twin primes in B which are larger than A. That is because T(A) contains all possible potential twin primes in A. Any other integers in A that are not potential twin primes are excluded from T(A) because (1) they are composites or (2) they are adjacent to a composite under the core seed primes of A. Consequently, when the additional core seed prime of B is added to each integer's modular signature in A, those that are composites in A cannot become potential twin primes in B.

| | A | B | C | Pa | Pb | Ps | Pc | Avg T(A) in C | Avg T(B) in C | Pz Smallest Non-Core Seed Prime of C | (Pz)^2 |
|---|---|---|---|---|---|---|---|---|---|---|---|
| | Primorials | | | Core Seed Primes Used for Product Factors | | | | Average Potential Twin Primes | | | |
| 1 | 210 | 2,310 | 47# | 7 | 11 | 13 | 47 | 5 | 59 | 53 | 2,809 |
| 2 | 2,310 | 30,030 | 173# | 11 | 13 | 17 | 173 | 35 | 456 | 179 | 32,041 |
| 3 | 30,030 | 510,510 | 709# | 13 | 17 | 19 | 709 | 286 | 4,868 | 719 | 516,961 |
| 4 | 510,510 | 9,699,690 | 3,109# | 17 | 19 | 23 | 3,109 | 3,272 | 62,162 | 3,119 | 9,728,161 |
| 5 | 9,699,690 | 223,092,870 | 14,929# | 19 | 23 | 29 | 14,929 | 43,632 | 1,003,543 | 14,939 | 223,173,721 |
| 6 | 223,092,870 | 6,469,693,230 | 80,429# | 23 | 29 | 31 | 80,429 | 727,428 | 21,095,426 | 80,447 | 6,471,719,809 |
| 7 | 6,469,693,230 | 200,560,490,130 | 447,829# | 29 | 31 | 37 | 447,829 | 15,900,087 | 492,902,698 | 447,841 | 200,561,561,281 |
| 8 | 200,560,490,130 | 7,420,738,134,810 | 2,724,079# | 31 | 37 | 41 | 2,724,079 | 380,157,930 | 14,065,843,393 | 2,724,109 | 7,420,769,843,881 |

Table 20

However, because T($A_C$) and T($B_C$) are averages of the distribution of potential twin primes in each cycle of A in C and B in C, it is possible that the number of potential twin primes in one or more cycles of A in C or B in C could equal 0 for very large primorials, even with the limit on the largest prime factor of C. Since under Theorem 7, T($B_C$) – T($A_C$) > 5, there must be one or more cycles of B in C in which the number of potential twin primes would be greater than 5 to offset any cycles in which the number of potential twins equaled 0. Therefore, there will be at least one potential twin prime in B that is not in A and is larger than A. By Theorem 5, all of those potential twin primes of B in C will remain primes because they are less than the square of the smallest non-core seed prime of C, as shown in the last column of Table 20.

Building the Primorial Scaffold for Goldbach Solutions

Unlike the scaffold for the twin primes which used three primorials, the scaffold for the Goldbach solutions will only use two. The reason is simple: the proof for the twin primes requires showing that there is always a pair of twin primes <u>greater</u> than some number N (so that they are infinite in number), whereas the proof for the Goldbach conjecture requires showing that there is a pair of odd primes <u>less</u> than any even number that will sum to that even number.

Before building the scaffold for the Goldbach solutions, Theorem 8 disposes of the case in which, for an even number E, there is an odd prime integer less than E/2 whose modular residues all differ from the corresponding modular residues of E under the smallest primorial greater than E. After that, we build the scaffold and prove that such a prime must exist.

**Theorem 8**: Let E equal a positive even integer > 4, and P1 equal an odd prime integer less than E/2. If the modular residues of E all differ from the corresponding modular residues of P1 under the seed primes for the smallest primorial, A, greater than E, then P2 = E – P1 will be a prime number, and P1 and P2 will form a Goldbach solution for E.

Proof.
1. Let P = the set of seed primes, {2, 3, … $P_N$}, of the smallest primorial, A, greater than E and let set P have C elements, where $P_N$ is the largest prime < $\sqrt{A}$.
2. Let $r_i$ = the modular residues of E for each $P_i \in P$, i = 1, 2, … C.
3. Let $s_i$ = the modular residues of P1 for each $P_i \in P$, with $r_i \neq s_i$ for each i.
4. Let $t_i$ equal the modular residue of P2 = E – P1 for each $P_i \in P$. Then $t_i = r_i - s_i \neq 0$ for each i.
5. Given that the elements of P are all prime numbers, then under the <u>Chinese Remainder Theorem</u> a unique integer solution (mod $\prod P$) exists to the set of equations:

$X = t_1 \bmod 2$

$X = t_2 \bmod 3$

- 
- 
- 

$X = t_C \bmod P_N$

6. Since X = P2 < E, it has a unique modular signature under the seed primes in P. Further, since each modular residue of P2 is not equal to zero for each $P_i \in P$, by <u>Theorem 2</u> P2 is a prime number. Further, since P2 = E – P1 is a prime number, it is co-prime to all prime numbers larger than $P_N$, so P1 and E cannot be in the same congruence class for any prime number greater than $P_N$ (or P2 would have residue = 0 for that prime number).
7. Since E = P1 + P2, and P1 and P2 are odd prime numbers, P1 and P2 will form a Goldbach solution for E.

Consider now what happens if we stack two primorials, A and B, where the largest core seed prime of B, Pb, is the largest prime $< \sqrt{A}$, and Pz is the smallest non-core seed prime of B. Table 21 below shows that every potential Goldbach solution of A in B will be less than the square of Pz.

| | Primorials | | | Largest Core Seed Prime | | Potential Goldbach Solutions of Primorial A | | | Pz Smallest Non-Core Seed Prime of B | |
|---|---|---|---|---|---|---|---|---|---|---|
| | A | B | Pa | Pb | T(A) in A | Product Factor | Avg T(A) in B | | (Pz)^2 |
| 1 | 210 | 13# | 7 | 13 | 15 | 0.692308 | 10 | 17 | 289 |
| 2 | 2,310 | 47# | 11 | 47 | 135 | 0.436373 | 59 | 53 | 2,809 |
| 3 | 30,030 | 173# | 13 | 173 | 1,485 | 0.307356 | 456 | 179 | 32,041 |
| 4 | 510,510 | 709# | 17 | 709 | 22,275 | 0.218553 | 4,868 | 719 | 516,961 |
| 5 | 9,699,690 | 3,109# | 19 | 3,109 | 378,675 | 0.164156 | 62,162 | 3,119 | 9,728,161 |
| 6 | 223,092,870 | 14,929# | 23 | 14,929 | 7,952,175 | 0.126197 | 1,003,543 | 14,939 | 223,173,721 |
| 7 | 6,469,693,230 | 80,429# | 29 | 80,429 | 214,708,725 | 0.098251 | 21,095,426 | 80,447 | 6,471,719,809 |
| 8 | 200,560,490,130 | 447,829# | 31 | 447,829 | 6,226,553,025 | 0.079161 | 492,902,698 | 447,841 | 200,561,561,281 |
| 9 | 7,420,738,134,810 | 2,724,079# | 37 | 2,724,079 | 217,929,355,875 | 0.064543 | 14,065,843,393 | 2,724,109 | 7,420,769,843,881 |

Table 21

Table 21 can be extended infinitely, of course, since there are an infinite number of primorials. Since all Goldbach solutions for E ≤ 210 have been accounted for in Table 10, consider that every even number E > 210 can be slotted between two successive primorials as shown in column 2 of Table 21 (extended infinitely). For example, if E = 250,000, then A = 30,030 (the largest primorial < 250,000) and B equals primorial 173#. The smallest non-core seed prime of B equals 179, and $179^2$ > 30,030. Then by Theorem 5, all 456 potential Goldbach solutions of A in B will remain primes (less 1) and there will exist at least one Goldbach solution for E, with P1 < 30,030.

**Theorem 9**: Let E equal any positive even integer > 30. Then there exists an odd prime number, P1, less than E whose modular residues all differ from the corresponding modular residues in the modular signature of E under the seed primes for E, where those seed primes are the prime factors of a primorial greater than E.

Proof.
1. Assume P1 does not exist.
2. Let A equal the largest primorial less than or equal to E.
3. Let $P_A$ equal the largest prime factor of A.
4. Let $P_B$ equal the largest prime $< \sqrt{A}$.
5. Let B equal the primorial $P_B$#. $P_B$ will be the largest core seed prime of B.
6. Let $P_Z$ equal the smallest prime > $P_B$. <u>$P_Z$ will be the smallest non-core seed prime of B.</u>

7. Let P = the set of seed primes of B and let P have N elements. Let $r_i$ = the modular residues of E in its modular signature under primorial B, (i.e., for each $P_i \in P$, i = 1, 2, … N, $r_i$ = E Mod $P_i$ in the modular signature of E).
8. Using Equation (4), let T(B) = (3-2)(5-2)…($P_A$ – 2)…($P_B$ – 2)

    T(B) will equal the number of combinations of the core seed primes > 2 of primorial B which do not contain 0 or the corresponding modular residue of E for a core seed prime of B, starting with prime 3 (because there will be only one residue available for 2 for all odd numbers). T(B) is the number of potential Goldbach solutions in primorial B.

9. Let GSet-B = the set of the potential Goldbach solutions for B, all of which will be odd potential prime numbers. GSet-B ≠ { } since every factor of T(B) > 0 and T(B) > 1.
10. Let GSet-A = the set of the potential Goldbach solutions for A in GSet-B. That is, GSet-A contains all potential Goldbach solutions of A that have not been converted to composites by the core seed primes of B.
11. By <u>Lemma 6.3</u>, the size of GSet-A ≥ 10 since the average number of potential Goldbach solutions in each cycle of A in B is ≥ 10.
12. Since $P_Z$ = the smallest prime > $P_B$, it is the smallest non-core seed prime of B.
13. By <u>Theorem 5</u>, every potential Goldbach solution of A in GSet-A that is less than $(P_Z)^2$ will remain a prime.
14. Since $(P_Z)^2$ > A, and GSet-A ≠ { }, there exists at least one prime P1 in GSet-A which is less than E and whose modular residues all differ from the corresponding modular residues of E under the seed primes for E.

## Part 5: Proof of the Twin Primes Conjecture

To prove the twin primes conjecture, we will need to do the following: (1) Assume there is a largest prime number, N, that is the last twin prime and (2) create a scaffold of three primorials greater than N, such that the middle primorial is large enough to contain potential twin primes greater than N and the square of the smallest non-core seed prime of the largest primorial is greater than the second primorial.

Proof.
1. Assume N > 31 since there are twin primes ≤ 31.
2. Let $T_W$ = the set of all twin primes and <u>assume that $T_W$ is a finite set</u>. Let $K_N$ = the number of elements in $T_W$. Then $K_N$ is the maximum number of twin primes among all positive integers.
3. Assume N is the largest odd prime which has a twin prime equal to N-2. That is, any other prime greater than N does not have a twin prime.
4. Let A = the smallest primorial > N, and let $P_A$ = the largest prime factor of A. That is,

    A = 2*3*5*…$P_A$.

5. Let B = the smallest primorial > A, and let $P_B$ = the largest prime factor (and core seed prime) of B. $P_B$ is also the smallest prime > $P_A$. That is,

    B = 2*3*5*…$P_A$*$P_B$.

6. Let $P_S$ = the smallest prime > $P_B$. $P_S$ will be the smallest non-core seed prime of B.
7. Let $P_C$ equal the largest prime number < $\sqrt{B}$. $P_C$ will be the largest non-core seed prime of B.
8. Let C = the primorial whose largest factor is $P_C$. That is,

    C = 2*3*5*…$P_A$*$P_B$*$P_S$…$P_C$

9. Using Equation (4), let T(C) = the number of potential twin primes of C under the core seed primes of C.

$$T(C) = (3-2)(5-2)(7-2)\ldots(P_C - 2)$$

10. Let $P_{NC}$ = the smallest non-core seed prime of C.
11. Let Q = the set of prime factors of C that are greater than $P_A$, and let Z = n(Q).
12. Using Equation (4), let $T(A_A)$ = the number of potential twin primes of A under the core seed primes of A.

$$T(A_A) = (3-2)(5-2)(7-2)\ldots(P_A - 2)$$

12. By <u>Theorem 6</u>, $T(A_C)$ = the average number of potential twin primes of A under the core seed primes of C,

$$T(A_C) = T(A_A) * \prod_{r=1}^{Z}((Pr - 2)/Pr) \text{ where } Pr \in Q$$
$$T(A_C) = T(C) / \prod_{r=1}^{Z}(Pr)$$

13. Using Equation (4), let $T(B_B)$ = the number of potential twin primes of B in B

$$T(B_B) = (3-2)(5-2)(7-2)\ldots(P_A - 2)(P_B - 2)$$

14. By <u>Theorem 6</u>, $T(B_C)$ = the average number of potential twin primes of B under the core seed primes of C,

$$T(B_C) = T(B_B) * \prod_{r=2}^{Z}((Pr - 2)/Pr) \text{ where } Pr \in Q$$
$$T(B_C) = T(C) / \prod_{r=2}^{Z}(Pr)$$

15. $T(B_C) - T(A_C)$ is the difference in the average number of potential twin primes between B and A in primorial C.

    Because of the assumption that $K_N$ is the maximum number of twin primes among all positive integers, all potential twin primes of A in C and B in C that are greater than N must be converted to composites by zero residues of the core seed primes of C. Then $T(B_C) - T(A_C) = 0$ since there will be no potential twin primes in B in C that are greater than A in C.

16. By <u>Theorem 7</u>, $T(B_C) - T(A_C) > 5$
17. Since $T(B_C) - T(A_C) > 5$, there are potential twin primes of B in C which are greater than A > N.
18. Since $P_{NC} > \sqrt{B}$, $(P_{NC})^2 > B$.
19. Since $(P_{NC})^2 > B$, by <u>Theorem 5</u> no potential twin prime of B in C will be converted to a composite, including every potential twin prime in B > A > N generated by the equation for $T(B_C)$.
20. Since every potential twin prime of B in C will remain a prime, and there are potential twin primes in B > A > N, then there exist twin primes that are greater than any element in $T_W$. Since N can be any odd prime number, the assumption is false that $T_W$ is a finite set. Therefore, there are an infinite number of twin primes.

## Part 6: Proof of the Goldbach Conjecture

There are three cases to the proof.

1. Case 1 (the trivial case): For E > 4 and E = 2N, where N = an odd prime number. If N is an odd prime number, then N + N = E, satisfying the requirement. Example: E = 6 = 2*3. Therefore, 3 + 3 = 6 satisfies the Goldbach requirement.

2. For E > 4 and E = 2N, where N ≠ an odd prime number.
   - Let ErP = the set of even numbers, $E_i$, for which $E_i$ = the product of consecutive prime numbers ≥ 2. That is, ErP is the set of primorials. Since there are an infinite number of prime numbers, ErP is an infinite ordered set.
   - Let $E_{MAX}$ = the smallest primorial ≥ E.
   - Let $P_{SP}$ = the set of <u>seed primes</u> for $E_{MAX}$, consisting of all prime numbers ≤ $\sqrt{Emax}$. Accordingly, each prime number less than or equal to E will have a unique modular signature under $P_{SP}$, the set of seed primes of $E_{MAX}$.
   - Let M = the number of elements in $P_{SP}$ and let $P_M$ = the largest prime in $P_{SP}$.

   Case (2a): Let P1 = a <u>seed prime</u> that satisfies the condition that no modular residue in the modular signature of P1 matches the corresponding modular residue in the modular signature of E for all $P_r \in P_{SP}$.

   Then by <u>Theorem 8</u>, P2 = E – P1 will be a prime number, and P1 and P2 will form a Goldbach solution for E.

   Case (2b): Assume no seed prime $P_r \in P_{SP}$ satisfies the condition that no modular residue in the modular signature of $P_r$ matches the corresponding modular residue in the modular signature of E for all $P_r \in P_{SP}$, where $P_r \leq P_N$. That is, for every $P_r \in P_{SP}$, one of its modular residues matches the corresponding modular residue of E. Therefore, assume no seed prime can be part of a Goldbach solution for E.

   1. Assume E = the smallest even number for which there is no matching pair of primes, P1 and P2, such that P1 + P2 = E.
   2. E must be greater than 30 since there are Goldbach solutions for E <= 30.
   3. Since there is no pair of prime numbers, P1 and P2, which satisfy the Goldbach conjecture for E, then there must be no prime number less than E, all of whose modular residues differ from the corresponding modular residues in the modular signature of E under the seed primes for E.
   4. However, that contradicts <u>Theorem 9</u>, so there exists P1 < E, where every modular residue in the modular signature of P1 will differ from the corresponding modular residue in the modular signature of E for every seed prime of E.
   5. Then P2 = E – P1 will have modular residues which all differ from those in E in its modular signature, and by <u>Theorem 2</u>, P2 will be a prime number since those modular residues are all non-zero. Then P1 and P2 are prime numbers whose sum is E and they form a Goldbach solution for E, in contradiction to the assumption. Then there is no E for which there is no Goldbach solution (i.e. for every even number, there is at least one pair of odd prime numbers less than E whose sum equals E).25